\documentstyle[12pt,leqno]{article}
\title{\ }
\author{\ }
\begin{document}
\newtheorem{theorem}{Theorem}[section]
\newtheorem{proposition}{Proposition}[section]
\newtheorem{lemma}{Lemma}[section]
\newtheorem{definition}{Definition}[section]
\newtheorem{corollaire}{Corollary}[section]
\newtheorem{remarque}{Remark}[section]
\def\Lim{\displaystyle\lim}
\def\Sup{\displaystyle\sup}
\def\Inf{\displaystyle\inf}
\def\Max{\displaystyle\max}
\def\Min{\displaystyle\min}
\def\Sum{\displaystyle\sum}
\def\Frac{\displaystyle\frac}
\def\Int{\displaystyle\int}
\def\n{|\kern -.05cm{|}\kern -.05cm{|}}
\def\Bigcap{\displaystyle\bigcap}
\def\E{{\cal E}}
\def\R{{\bf \hbox{\sc I\hskip -2pt R}}} 
\def\N{{\bf \hbox{\sc I\hskip -2pt N}}} 
\def\Z{{\bf Z}} 
\def\Q{{\bf \hbox{\sc I\hskip -7pt Q}}} 
\def\C{{\bf \hbox{\sc I\hskip -7pt C}}} 

{\begin{center} {\ } \vskip 2.1cm {\large\bf Indirect boundary stabilization for weakly  coupled degenerate
wave equations under fractional damping
 \\ } \
\\ \ \\
\begin{tabular}{c}
\sc Rachid Benzaid and Abbes Benaissa,
\\ \\
\small  Laboratory of Analysis and Control of PDEs, \qquad \\
\small Djillali Liabes University,\qquad \\
\small P. O. Box 89, Sidi Bel Abbes 22000, ALGERIA.\qquad \\
\small E-mails: n.mouslim@yahoo.fr\\
\small rachidbenzaid995@gmail.com \\
\end{tabular}
\end{center}
\
\\
\begin{abstract}
{\small In this paper, we consider the well-posedness and stability of a one-dimensional system
of degenerate wave equations coupled via zero order terms with one boundary fractional damping acting on one end only. We prove optimal
polynomial energy decay rate of order $1/t^{(3-\tau)}$. The method is based on
the frequency domain approach combined with multiplier technique.}
\end{abstract}


{\it ${\cal AMS}$ Classification:} 35B40,35L80, 74D05, 93D15.
\begin{sloppypar}
\section{Introduction}
In this paper, we investigate the existence and energy decay rate of a system of
coupled degenerate wave equations with only one fractional boundary damping.
This system defined on $(0,1)\times (0,+\infty )$ takes the
following form

\begin{equation}
\left\{
\matrix{
u_{tt}(x,t)-(a(x)u_{x})_{x}(x,t)+\alpha v=0 \hfill & \hbox{ in }(0,1)\times
(0,+\infty ), \hfill & \cr
v_{tt}(x,t)-(a(x)v_{x})_{x}(x,t)+\alpha u=0\hfill & \hbox{ in }(0,1)\times
(0,+\infty ),\hfill & \cr
\left\{\matrix{u(0, t)=0 \hfill & \hbox{ if } 0\leq m_a< 1 \hfill & \cr
(a u_{x})(0,t)=0 \hfill & \hbox{ if } 1\leq m_a< 2 \hfill & \cr}\right.\hfill  & \hbox{ in } (0, +\infty),\hfill&\cr
\left\{\matrix{v(0, t)=0 \hfill & \hbox{ if } 0\leq m_a< 1 \hfill & \cr
(a v_{x})(0,t)=0 \hfill & \hbox{ if } 1\leq m_a< 2 \hfill & \cr}\right.\hfill  & \hbox{ in } (0, +\infty),\hfill&\cr
v(1, t)=0 \hfill & \hbox{ for } t\in (0, +\infty),\hfill & \cr
\beta u(1,t)+ (a u_{x})(1, t)=- \varrho \partial_{t}^{\tau, \omega} u(1, t) \hfill & \hbox{ in } (0, +\infty),\hfill&\cr
u(x,0)=u_0(x), u_{t}(x,0)=u_1(x),
v(x,0)=v_0(x), v_{t}(x,0)=v_1(x)\hfill & \hbox{ for } x\in (0, 1),\hfill &\cr}
\right.
\label{1}
\end{equation}
where $a\in C([0, 1])\cap C^1(]0, 1])$ is positive on $]0, 1]$ but vanishes at zero, $\alpha$
denote the coupling parameter, which is assumed to be real and small enough, $\beta\geq$ 0 and $\varrho >0$.
The notation $\partial^{\tau,\omega}_{t}$ stands for
the generalized Caputo's fractional derivative of order $\tau, $  $(0 < \tau< 1) $, with respect to the
time variable (see {\bf\cite{choi}}). It is defined as follows
$$
\partial^{\tau,\omega}_{t}g(t)=\left\{
\matrix{
g_{t} & \hbox{for} \ \tau=1, \ \omega\geq 0,\hfill &\cr
\frac{1}{\Gamma(1-\tau)}\int_{0}^{t}(t-s)^{-\tau}e^{-\omega(t-s)}  \frac{dg}{ds}(s)ds\quad &
\hbox{for }\ 0<\tau <1,\quad\omega\geq 0. \ \hfill &\cr}
\right.
$$
The initial data $(u_{0},u_{1},v_{0},v_{1})$ belong to a suitable
function space.

Degenerate partial differential equations are encountered in the theory of boundary layers, in the theory of shells,
in the theory of diffusion processes, in particular in the theory of Brownian motion, in climate science,
in contact mechanics and in many other problems in physics and mechanics.
We find that the commun feature of these problems is the lose of its typical characteristics,
including ellipticity or hyperbolicity, which can have a substantial impact on how
solutions behave.

Degenerate equations are studied by posing two closely connected problems: 1) a demonstration of the solvability of,
say, boundary value problems taking into account changes in their formulation which are a consequence of the
degeneration of type; and 2) a determination of properties of the solutions which are analogous
to those of non-degenerate equations (smoothness, Harnack inequalities for elliptic and parabolic equations, etc.).

We review the related papers, regarding linear degenerate wave system, from a qualitative and quantitative
study. For a single degenerate wave equation, we beginning with the work treated in {\bf\cite{alabau1}}, for
$(x, t)\in (0, 1)\times (0, +\infty)$ where the goal was mainely on the equation
$$
u_{tt}(x,t)-(a(x) u_{x}(x,t))_x=0 \hbox{ in } (0, 1)\times (0, \infty),
$$
together with boundary linear damping of the form
$$
\left\{\matrix{\left\{\matrix{u(0, t)=0 \hfill & \hbox{ if } 0\leq m_a< 1 \hfill & \cr
(a u_{x})(0,t)=0 \hfill & \hbox{ if } 1\leq m_a< 2 \hfill & \cr}\right.\hfill  & \hbox{ in } (0, +\infty) ,\hfill&\cr
u_t(1, t)+u_x(1, t)+\beta u(1, t)=0 \hfill & \hbox{ in } (0, +\infty).\hfill&\cr
}\right.
$$
where $\beta>0$ is the given constant. $m_a=\sup_{0< x\leq 1}\Frac{x|a'(x)|}{a(x)}<2$ is the measurement
of the degree of the degeneracy.
Thanks to the energy multiplier method, it is proved that the total energy of the whole
system decays exponentially.

Recently,  Benaissa and Aichi {\bf\cite{beai}} considered
the scalar degenerate wave equation under the following  boundary fractional damping
$$
\left\{\matrix{\left\{\matrix{u(0, t)=0 \hfill & \hbox{ if } 0\leq m_a< 1 \hfill & \cr
(a u_{x})(0,t)=0 \hfill & \hbox{ if } 1\leq m_a< 2 \hfill & \cr}\right.\hfill  & \hbox{ in } (0, +\infty) ,\hfill&\cr
(au_{x})(1, t)+\varrho\partial_{t}^{\tau, \omega}u(1, t)+\beta u(1,t)=0 \hfill & \hbox{ in } (0, +\infty).\hfill&\cr
}\right.
$$
They obtained an optimal polynomial stability of the solutions by using a frequency domain approach combining
with a multiplier method.

Next, in a recent paper of Liu and Rao {\bf\cite{lira}} general systems of coupled second
order evolution equations have been studied. The system is described
$$
\left\{\matrix{u_{tt}-b\Delta u+\alpha y=0 \hfill & \hbox{ on } \Omega,\hfill \cr
y_{tt}-\Delta u+\alpha u=0 \hfill & \hbox{ on } \Omega,\hfill \cr
u=0 \hfill  &\hbox{ on } \Gamma_D, \hfill \cr
b\partial_{\nu} u+\gamma u+u_t=0\hfill &\hbox{ on }\Gamma_N,\hfill\cr
y=0 \hfill&\hbox{ in } \Gamma,\hfill\cr}\right.
$$
where $\Omega\subset \R ^n$ is a bounded domain with smooth boundary $\Gamma$ of class $C^2$
such that $\Gamma=\Gamma_D\cup \Gamma_N$ and $\Gamma_D\cap \Gamma_N=\emptyset$.
They established, by the frequency domain approach, polynomial decay
rate of order $\frac{\ln t}{t}$ for smooth initial data, while waves propagate with equal speeds.
Moreover, while waves propagate with different speeds, i.e. the case $b\not= 1$, they
proved that the energy decays at a rate which depends on the arithmetic property
of the ratio of the wave speeds $b$.

Very recently, Wehbe and Koumaiha {\bf\cite{weko}} considered a one-dimensional setting of a system of wave equation coupled via zero order terms.
More precisely, they studied the stabilization
of the following system  of partially damped coupled wave equations propagating with equal speeds, described by
$$
\left\{
\matrix{
u_{tt}-u_{xx}+\alpha y=0 \hfill & \hbox{ in }(0,1)\times (0,+\infty ), \hfill & \cr
y_{tt}-y_{xx}+\alpha u=0\hfill & \hbox{ in }(0,1)\times (0,+\infty ),\hfill & \cr
u(0, t)=y(0, t)=y(1, t)=0 \hfill  & \hbox{ in } (0, +\infty) ,\hfill&\cr
u_{x}(1, t)+\gamma u_t(1, t)=0 \hfill & \hbox{ in } (0, +\infty),\hfill&\cr
u(x,0)=u_0(x), u_{t}(x,0)=u_1(x),
y(x,0)=y_0(x), y_{t}(x,0)=y_1(x)\hfill & \hbox{ for } x\in (0, 1),\hfill &\cr}
\right.
$$
where $\gamma> 0$. They proved optimal polynomial energy decay rate of order $\frac{1}{t}$, by using a
frequency domain approach and Riesz basis property of the generalized eigenvector of the system.

In {\bf\cite{akil}}, Akil et al considered a one-dimensional coupled wave equations on its indirect boundary
stabilization defined by
$$
\left\{
\matrix{
u_{tt}(x,t)-u_{xx} (x,t)-d v_t(x,t)=0 \hfill &\hbox{ in }
(0, 1)\times (0,+\infty ), \hfill & \cr
v_{tt}(x,t)-v_{xx}+d u_t(x,t)=0\hfill & \hbox{ in }%
(0, 1)\times (0,+\infty ),\hfill & \cr
u(0,t)=v(0,t)=v(1, t)=0 \hfill & \hbox{ on } (0,+\infty ),\hfill & \cr
u_{x}(1,t)+\varrho \partial _{t}^{\tau ,\omega }u(1,t)=0
\hfill & \hbox{ on } (0,+\infty ),\hfill & \cr
u(x,0)=u_{0}(x),\ u_{t}(x,0)=u_{1}(x)\hfill & \hbox{ on }
(0,1),\hfill & \cr
v(x,0)=v_{0}(x),\ v_{t}(x,0)=v_{1}(x)\hfill & \hbox{ on }(0,1).\hfill & \cr}
\right.
$$
They established a polynomial energy decay rate of type $t^{-s(\tau)}$, such that

\noindent
i) If $d\not= k\pi$, then $s(\tau)=\frac{2}{1-\tau}$.

\noindent
ii)If $d= k\pi$, then $s(\tau)=\frac{2}{5-\tau}$.\\

In {\bf\cite{1}}, kerdache et al investigate the decay rate of the energy
of the coupled wave equations with a two boundary fractional dampings,
that is,
$$
\left\{
\matrix{
u_{tt}(x,t)-u_{xx} (x,t)+\alpha (u-v)=0 \hfill &\hbox{ in }
(0, 1)\times (0,+\infty ), \hfill & \cr
v_{tt}(x,t)-v_{xx}+\alpha (v-u)=0\hfill & \hbox{ in }%
(0, 1)\times (0,+\infty ),\hfill & \cr
u(0,t)=v(0,t)=0 \hfill & \hbox{ on } (0,+\infty ),\hfill & \cr
u_{x}(1,t)+\varrho \partial _{t}^{\tau ,\omega }u(1,t)=0
\hfill & \hbox{ on } (0,+\infty ),\hfill & \cr
v_{x}(1,t)+\tilde{\varrho}\partial _{t}^{\tau ,\omega
}u(1,t)=0\hfill & \hbox{ on }%
(0,1)\times (0,+\infty ),\hfill & \cr
u(x,0)=u_{0}(x),\ u_{t}(x,0)=u_{1}(x)\hfill & \hbox{ on }
(0,1),\hfill & \cr
v(x,0)=v_{0}(x),\ v_{t}(x,0)=v_{1}(x)\hfill & \hbox{ on }(0,1).\hfill & \cr}
\right.
$$
Using semigroup theory, they proved an optimal polynomial type decay
rate.

Motivated by the works {\bf\cite{lira}}, {\bf\cite{beai}} and {\bf\cite{weko}} we wonder what the asymptotic behavior of the coupled degenerate wave equations
would be, considering a boundary fractional damping acting only on one equation.

This paper is divided into four sections. In section 2,
we introduce the appropriate functional spaces that are naturally associated with degenerate problems and
preliminary result used throughout the paper.
Section 3 is devoted to the proof of the well-posedness and strong asymptotic of the considered system.
In Section 4 we establish an optimal polynomial decay of type $t^{-\frac{2}{3-\tau}}$ for smooth initial
data, by the frequency domain method.

\section{Preliminary results}
Let $a\in C([0, 1]\cap C^1(]0, 1])$ be a function satisfying the following assumptions:
\begin{equation}
\left\{\matrix{(i)\hfill & a(x)> 0\ \forall x\in ]0, 1], a(0)=0, \hfill &\cr
(ii) \hfill &m_a=\Sup_{0< x\leq 1}\Frac{x|a'(x)|}{a(x)}<2, \hbox{ and } \hfill &\cr
(iii) \hfill & a\in C^{[m_a]}([0, 1]),\hfill &\cr}\right.
\label{ee211}
\end{equation}
where $[\cdot]$ stands for the integer part.

When $m_a>1$, we suppose $\beta>0$ because if $\beta=0$ and the feedback law only depends on velocities, we
may encounter the situation where the closed-loop system is not well-posed in terms
of the semigroups in the Hilbert space.\\

\noindent
{\bf Examples:} 1) Let $\varpi\in (0, 2)$ be given. Define
$$
a(x)=x^{\varpi} \quad \forall x\in [0, 1].
$$
satifies (\ref{ee211}).

\noindent
2) Let $\varpi\in [0, 2)$ be given and let $\theta\in (0, 1-\varpi/2)$. The function
$$
a(x)=x^{\varpi}(1+\cos^2(\ln x^{\theta})) \quad \forall x\in [0, 1]
$$
satifies (\ref{ee211}).

Now, we introduce, as in {\bf\cite{canna.3}}, {\bf\cite{foto}} or {\bf\cite{alabau1}}, the following weighted spaces:
$$
H_{a}^{1}(0,1)=\left\{u \hbox{ is locally absolutely continuous in } (0, 1]: \sqrt{a(x)}u_x \in L^2(0,1)\right\}.
$$
It is easy to see that $H_a^1(0, 1)$ is a Hilbert space with the scalar product
$$
(u, v)_{H_a^1(0, 1)}=\Int_{0}^{1}(a(x) u'(x)\overline{v'(x)}+u(x)\overline{v(x)})\, dx\quad \forall u, v\in H_a^1(0, 1)
$$
and associated norm
$$
\|u\|_{H_a^1(0, 1)}=\left\{\Int_{0}^{1}(a(x) |u'(x)|^2+|u(x)|^2)\, dx\right\}^{1/2}\quad \forall u\in H_a^1(0, 1).
$$
Next, we define
$$
H_{a}^2(0, 1)=\{u\in H_{a}^1(0, 1):\ au'\in H^1 (0, 1)\},
$$
where $H^1 (0, 1)$ denotes the classical Sobolev space.

In order to express the boundary conditions of the first component of the solution of
(\ref{1}) in the functional setting, we define the spaces $H_{0, a}^{1}(0,1)$ and $W_{a}^1(0, 1)$ depending on the value of $m_a$,
as follows:
\begin{itemize}
\item[(i)] For $0\leq m_a< 1$, we define
$$
\left\{\matrix{H_{0, a}^{1}(0,1)=\left\{u\in H_{a}^{1}(0, 1) /\ u(0)=u(1)=0\right\}, \hfill & \cr
W_{a}^{1}(0,1)=\left\{u\in H_{a}^{1}(0, 1) /\ u(0)=0\right\}. \hfill & \cr}\right.
$$
\item[(ii)] For $1\leq m_a< 2$, we define
$$
\left\{\matrix{H_{0, a}^{1}(0,1)=\left\{u\in H_{a}^{1}(0, 1) /\ u(1)=0\right\}, \hfill & \cr
W_{a}^{1}(0,1)=H_a^{1}(0,1). \hfill & \cr}\right.
$$
\end{itemize}
It is easy to see that $H_a^1(0, 1)$ when $\beta> 0$ is a Hilbert space with the scalar product
$$
(u, v)_{H_a^1(0, 1)}=\Int_{0}^{1}a(x) u'(x)\overline{v'(x)}\, dx +\beta u(1)\overline{v(1)}.
$$
Let us also set
$$
|u|_{*}=\left(\Int_{0}^{1}a(x) |u'(x)|^2\, dx\right)^{1/2}\quad \forall u\in H_a^1(0, 1).
$$
Actually, $|\cdot|_{*}$ is an equivalent norm on the closed subspaces $H_{0, a}^1(0, 1)$ and $W_{a}^{1}(0,1)$ to the
norm of $H_{a}^1(0, 1)$ when $m_a\in [0, 1[$.
This fact is a simple consequence of the following version of Poincar\'e's inequality.
\begin{proposition} Assume (\ref{ee211}) with $m_a\in [0, 1)$. Then there is a positive constant $C_*=C(a)$ such that
\begin{equation}
\|u\|_{L^{2}(0, 1)}^2\leq C_* |u|_{1, a}^2\quad \forall u\in H_{0, a}^1(0, 1).
\label{eh10}
\end{equation}
\end{proposition}
{\bf Proof.}
Let $u\in H_{0, a}^1(0, 1)$. For any $x \in ]0, 1]$ we have that
$$
|u(x)|=\left|\Int_{0}^{x}u'(s)\, ds\right|\leq |u|_{1, a}\left\{\Int_{0}^{1}\Frac{1}{a(s)}\, ds\right\}^{1/2}.
$$
Therefore
$$
\Int_{0}^{1}|u(x)|^2\, dx\leq |u|_{1, a}^2\left\{\Int_{0}^{1}\Frac{1}{a(s)}\, ds\right\}.
$$
Now, we state two propositions that will be needed later (see {\bf\cite{canna.3}}, {\bf\cite{foto}} and {\bf\cite{alabau1}}).
\begin{proposition} Assume (\ref{ee211}). Then the following properties hold.
\begin{itemize}
\item[(i)] For every $u\in H_{a}^1(0, 1)$
\begin{equation}
\Lim_{x\rightarrow 0}x u^{2}(x)=0.
\label{eh100}
\end{equation}
\item[(ii)] For every $u\in H_{a}^2(0, 1)$
\begin{equation}
\Lim_{x\rightarrow 0}x a(x)u'(x)^{2}=0.
\label{eh101}
\end{equation}
\item[(iii)] For every $u\in H_{a}^2(0, 1)$
\begin{equation}
\Lim_{x\rightarrow 0}x a(x)u(x) u'(x)=0.
\label{eh10199}
\end{equation}
\end{itemize}
\label{propo4}
\end{proposition}
\begin{proposition}
$H_{a}^1(0, 1)\hookrightarrow L^{2}(0, 1)$ with compact embedding.
\end{proposition}

\section{Well-posedness and strong stability}
\subsection{Augmented model}
In this section we reformulate $(P)$ into an augmented system. For that, we need the following proposition.
\begin{proposition}[see {\bf\cite{1}}]
Let $\vartheta$ be the function:
\begin{equation}
\vartheta(\varsigma)=|\varsigma|^{(2\tau-1)/2},\quad -\infty<\varsigma<+\infty,\
0<\tau<1. \label{e1}
\end{equation}
Then the relationship between the `input' U and the `output' O of
the system
\begin{equation}
\partial_t\varphi(\varsigma, t)+(\varsigma^{2}+\omega)\varphi(\varsigma, t) -U(t)\vartheta(\varsigma)=0,\quad -\infty<\varsigma<+\infty,\omega\geq 0, t> 0,
\label{e6}
\end{equation}
\begin{equation}
\varphi(\varsigma, 0)=0, \label{e7}
\end{equation}
\begin{equation}
O(t)=(\pi)^{-1}\sin(\tau\pi)\Int_{-\infty}^{+\infty}\vartheta(\varsigma)\varphi(\varsigma, t)\, d\varsigma,
\label{e8}
\end{equation}
where $U\in C^0([0,+\infty))$, is given by
\begin{equation}
O=I^{1-\tau, \omega}U, \label{e9}
\end{equation}
where
$$
[I^{\tau,
\omega}f](t)=\Frac{1}{\Gamma(\tau)}\Int_{0}^{t}(t-s)^{\tau-1}e^{-\omega(t-s)}f(s)\, ds.
$$
\label{th2}
\end{proposition}
\begin{lemma}[see {\bf\cite{1}}]
If $\lambda\in D_{\omega}=\C\backslash]-\infty, -\omega]$ then
$$
\Int_{-\infty}^{+\infty}\Frac{\vartheta^2(\varsigma)}{\lambda+\omega+\varsigma^2}\, d\varsigma
=\Frac{\pi}{\sin\tau\pi}(\lambda+\omega)^{\tau-1}.
$$
\label{achour}
\end{lemma}

We are now in a position to reformulate system $(P)$. Indeed, by using Proposition \ref{th2}, system
$(P)$ may be recast into the augmented model:
$$
\left\{
\matrix{
u_{tt}(x,t)-(a(x)u_{x})_{x}(x,t)+\alpha v=0 \hfill & \hbox{ in }(0,1)\times
(0,+\infty ), \hfill & \cr
v_{tt}(x,t)-(a(x)v_{x})_{x}(x,t)+\alpha u=0\hfill & \hbox{ in }(0,1)\times
(0,+\infty ),\hfill & \cr
\varphi_t(\varsigma, t)+(\varsigma^{2}+\omega)\varphi(\varsigma, t) -u_{t}(1, t)\vartheta(\varsigma)=0, \hfill & -\infty <\varsigma<+\infty,\hfill \omega\geq 0,\hfill t>0,\hfill&\cr
\left\{\matrix{u(0, t)=0 \hfill & \hbox{ if } 0\leq m_a< 1 \hfill & \cr
(a u_{x})(0,t)=0 \hfill & \hbox{ if } 1\leq m_a< 2 \hfill & \cr}\right.\hfill  & \hbox{ in } (0, +\infty) ,\hfill&\cr
\left\{\matrix{v(0, t)=0 \hfill & \hbox{ if } 0\leq m_a< 1 \hfill & \cr
(a v_{x})(0,t)=0 \hfill & \hbox{ if } 1\leq m_a< 2 \hfill & \cr}\right.\hfill  & \hbox{ in } (0, +\infty) ,\hfill&\cr
v(1, t)=0 \hfill & \hbox{ for } t\in (0, +\infty),\hfill & \cr
\beta u(1,t)+ (a u_{x})(1, t)=-\zeta \Int_{-\infty}^{+\infty}\vartheta(\varsigma)\varphi(\varsigma, t)\, d\varsigma, \hfill & \zeta=\varrho(\pi)^{-1}\sin(\tau\pi),\hfill \cr
u(x,0)=u_0(x), u_{t}(x,0)=u_1(x),
v(x,0)=v_0(x), v_{t}(x,0)=v_1(x)\hfill & \hbox{ for } x\in (0, 1).\hfill &\cr}
\right.
\leqno{(P')}
$$
We define the energy associated to the solution of the problem $(P')$ by the following formula:
\begin{equation}
\matrix{{\cal E}(t)=\Frac{1}{2}\int_{0}^{1}(|u_t|^2+a(x) |u_x|^2)dx+\Frac{1}{2}\int_{0}^{1}(|v_t|^2+a(x) |v_x|^2)dx\hfill &\cr
\Frac{1}{2}\alpha\int_{0}^{1}(u {\overline v}+v{\overline u})dx
+\Frac{\beta}{2} |u(1,t)|^{2}+\Frac{\zeta}{2}\int_{-\infty}^{+\infty}|\varphi(\varsigma, t)|^{2}\, d\varsigma.\hfill & \cr}
\label{e10}
\end{equation}
\begin{lemma}
Let $(u, v, \varphi)$ be a regular solution of the problem $(P')$.
Then, the energy functional defined by (\ref{e10}) satisfies
\begin{equation}
{\cal E}'(t)=-\zeta \int_{-\infty}^{+\infty}(\varsigma^{2}+\omega)|\varphi(\varsigma, t)|^{2}\, d\varsigma\leq 0.
\label{e11}
\end{equation}
\label{lem1}
\end{lemma}
In this section, we give an existence and uniqueness result for problem $(P')$ using the semigroup theory.
Introducing the vector function $U=(u, {\tilde u}, v, {\tilde v}, \varphi)^{T}$, where ${\tilde u}=u_t$, ${\tilde v}=v_t$,
system $(P')$ can be treated as a Cauchy evolution problem
\begin{equation}
\left\{\matrix{\Theta'={\cal A}\Theta,\qquad\hbox{ for all } t> 0, \hfill & \cr
\Theta(0)=\Theta_0,\hfill & \cr}
\right.
\label{e14}
\end{equation}
where $\Theta_0=(u_0, u_1, v_0, v_1, \varphi_{0})^{T}$ and
$$
{\cal A}: D({\cal A})\subset {\cal H}\rightarrow {\cal H}
$$
is the operator given by
\begin{equation}
{\cal A}\pmatrix{u\cr {\tilde u}\cr v\cr {\tilde v}\cr \varphi\cr}=
\pmatrix{{\tilde u}\cr (a(x)u_x)_x-\alpha v\cr {\tilde v}\cr (a(x)v_x)_x-\alpha u\cr
-(\varsigma^{2}+\omega)\varphi+{\tilde u}(1)\vartheta(\varsigma)\cr}.
\label{e1972}
\end{equation}
We introduce the following phase space (the energy space):
$$
{\cal H}=W_{a}^{1}(0,1)\times L^{2}(0,1)\times H_{0, a}^{1}(0,1)\times L^{2}(0,1)\times L^{2}(-\infty, +\infty),
$$
that is a Hilbert space with the following inner product
$$
\matrix{\langle U, {\tilde U} \rangle_{\cal H}=\Int_{0}^{1}a(x)u_{1x} {\overline u}_{2x} dx+\Int_{0}^{1}a(x)v_{1x} {\overline v}_{2x} dx
+\alpha\Int_{0}^{1} (u_1 {\overline v}_2+v_2 {\overline u}_1)\, dx\hfill & \cr
+\Int_{0}^{1}{\tilde u}_1 {\overline {\tilde u}}_2 dx +\Int_{0}^{1}{\tilde v}_1 {\overline {\tilde v}}_2 dx
+\zeta\Int_{-\infty}^{+\infty} \varphi_1  {\overline \varphi_2}  \, d\varsigma + \beta u_1(1){\overline u_2}(1),\hfill &\cr}
$$
for all $U= (u_1, {\tilde u}_1, v_1, {\tilde v}_1, \varphi_1)^{T}$ and $\tilde{U}=(u_2, {\tilde u}_2, v_2, {\tilde v}_2, \varphi_2)^{T}$.

The domain of ${\cal A}$ is
\begin{equation}
D({\cal A})=\left\{\matrix{(u, {\tilde u}, v,{\tilde v}, \varphi)^{T}\ \hbox{in } {\cal H}:
u\in H_{a}^{2}(0,1)\cap W_{a}^{1}(0,1), v\in H_{a}^{2}(0,1)\cap H_{0, a}^{1}(0,1),  \hfill \cr
{\tilde u}\in W_{a}^{1}(0,1), {\tilde v}\in H_{0, a}^{1}(0,1),
-(\varsigma^{2}+\omega)\varphi+{\tilde u}(1)\vartheta(\varsigma)\in L^{2}(-\infty, +\infty)  ,\hfill \cr
\beta u(1)+(a u_x)(1)
+\zeta\Int_{-\infty}^{+\infty}\vartheta(\varsigma)\varphi(\varsigma)\, d\varsigma=0, \hfill \cr
|\varsigma|\varphi\in L^{2}(-\infty, +\infty)\hfill \cr
}\right\}.
\label{e15}
\end{equation}
We have the following existence and uniqueness result.
\begin{theorem}[Existence and uniqueness]\
\begin{itemize}
\item[(1)] If $U_0\in D({\cal A})$, then system (\ref{e14}) has a unique strong solution with the following regularity,
$$
U\in C^0(\R_{+}, D({\cal A}))\cap C^1(\R_{+}, {\cal H}).
$$
\item[(2)] If $U_0\in {\cal H}$, then system (\ref{e14}) has a unique weak solution such that
$$
U\in C^0(\R_{+}, {\cal H}).
$$
\end{itemize}
\label{1th}
\end{theorem}
{\bf Proof.}\\
We use the semigroup approach. In what follows, we prove that ${\cal A}$ is monotone. For any $U\in D({\cal A})$
and using (\ref{e14}), (\ref{e11}) and the fact that
\begin{equation}
{\cal E}(t)=\Frac{1}{2}\|U\|_{\cal H}^{2},
\label{e16}
\end{equation}
we have
\begin{equation}
\Re\langle {\cal A} U,  U \rangle_{\cal H}=-\zeta \Int_{-\infty}^{+\infty}(\varsigma^{2}+\omega)|\varphi(\varsigma)|^{2}\, d\varsigma,
\label{e17}
\end{equation}
and therefore, ${\cal A}$ is dissipative.
Next, we prove that the operator $\lambda I-{\cal A}$ is surjective for $\lambda > 0$. More
precisely, given $G=(g_1, g_2, g_3, g_4, g_5)^{T}\in {\cal H}$, we will show that there is $U\in D({\cal A})$ such
that
\begin{equation}
(\lambda I-{\cal A})U=G.
\label{e18z}
\end{equation}
From Equation (\ref{e18z}), we get the following system of equations
\begin{equation}
\left\{\matrix{\lambda u-{\tilde u}=g_1,  \hfill & \cr
\lambda {\tilde u}-(a(x)u_x)_x+\alpha v=g_2, \hfill & \cr
\lambda v-{\tilde v}=g_3\hfill & \cr
\lambda {\tilde v}-(a(x)v_x)_x+\alpha u=g_4, \hfill & \cr
\lambda \varphi+(\varsigma^{2}+\omega)\varphi-{\tilde u}(1)\vartheta(\varsigma)=g_5.\hfill & \cr}\right.
\label{e18}
\end{equation}
Suppose $u, v$ are found  with the appropriate regularity. Then,
$(\ref{e18})_1$ and $(\ref{e18})_3$ yield
\begin{equation}
\left\{\matrix{{\tilde u}=\lambda u-g_1 \in W_{a}^{1}(0,1),  \hfill & \cr
{\tilde v}=\lambda v-g_3 \in H_{0, a}^{1}(0,1),  \hfill & \cr}\right.
\label{e19}
\end{equation}
By using $(\ref{e18})_2, (\ref{e18})_4$ and (\ref{e19}) it can easily be shown that $u, v$
satisfy
\begin{equation}
\left\{\matrix{\lambda^{2}u -(a(x)u_x)_x+\alpha v=g_2+\lambda g_1,\hfill & \cr
\lambda^{2}u -(a(x)u_x)_x+\alpha u=g_4+\lambda g_3.\hfill & \cr}\right.
\label{e21}
\end{equation}
Solving system (\ref{e21}) is equivalent to finding $(u, v)\in H_a^{2}(0,1)\cap W_{a}^{1}(0,1)\times
H_a^{2}(0,1)\cap H_{0, a}^{1}(0,1)$ such that
\begin{equation}
\left\{\matrix{\Int_{0}^{1}(\lambda^{2}u {\overline w}- (a(x)u_x)_x {\overline w})\, dx+\alpha \Int_{0}^{1}v{\overline w}\, dx
=\Int_{0}^{1}(g_2+\lambda g_1){\overline w}\, dx,\hfill & \cr
\Int_{0}^{1}(\lambda^{2}v {\overline y}- (a(x) v_x)_x {\overline y})\, dx
=\Int_{0}^{1}(g_4+\lambda g_3){\overline y}\, dx,\hfill & \cr}\right.
\label{e22}
\end{equation}
for all $(w, y)\in W_{a}^{1}(0,1)\times H_{0, a}^{1}(0,1)$. By using (\ref{e22}), the boundary
condition $(\ref{e15})_3$ and $(\ref{e18})_5$ the functions $u$ and $v$
satisfy the following system
\begin{equation}
\left\{\matrix{
\matrix{ \Int_{0}^{1}(\lambda^{2}u {\overline w}+ a(x)u_x
{\overline w_x})\, dx +\alpha \Int_{0}^{1}v{\overline w}\, dx+\beta u(1) {\overline w(1)}+{\tilde\zeta} {\tilde u}(1) {\overline w(1)}\, \hfill &\cr
\quad\ \qquad=\Int_{0}^{1}(g_2+\lambda g_1){\overline w}\, dx
-\zeta\Int_{-\infty}^{+\infty} \Frac{\vartheta(\varsigma)}{\varsigma^{2}+\omega+\lambda}g_5(\varsigma)\,  d\varsigma{\overline w(1)}, &\cr}\hfill & \cr
\matrix{ \Int_{0}^{1}(\lambda^{2}v {\overline y}+ a(x)v_x
{\overline y_x})\, dx+\alpha \Int_{0}^{1}u{\overline y}\, dx=\Int_{0}^{1}(g_4+\lambda g_3){\overline y}\, dx, &\cr}\hfill & \cr}\right.
\label{e23}
\end{equation}
where ${\tilde\zeta}=\zeta\Int_{-\infty}^{+\infty}
\Frac{\vartheta^{2}(\varsigma)}{\varsigma^{2}+\omega+\lambda}  \, d\varsigma$. Using again
$(\ref{e19})_1$, we deduce that
\begin{equation}
{\tilde u}(1)=\lambda u(1)-g_1(1). \label{e24}
\end{equation}
Inserting (\ref{e24}) into (\ref{e23}), we get
\begin{equation}
\left\{\matrix{
\matrix{ \Int_{0}^{1}(\lambda^{2}u {\overline w}+  a(x)u_x
{\overline w_x})\, dx+\alpha \Int_{0}^{1}v{\overline w}\, dx +(\lambda {\tilde\zeta}+\beta) u(1){\overline w(1)}\, \hfill &\cr
=\Int_{0}^{1}(g_2+\lambda g_1){\overline w}\, dx-\zeta \Int_{-\infty}^{+\infty}
\Frac{\vartheta(\varsigma)}{\varsigma^{2}+\omega+\lambda}g_5( \varsigma)\,  d\varsigma{\overline w(1)} +{\tilde\zeta} g_1(1) {\overline w(1)}, &\cr}\hfill &\cr
\Int_{0}^{1}(\lambda^{2}v {\overline y}+ a(x)v_x
{\overline y_x})\, dx+\alpha \Int_{0}^{1}u{\overline y}\, dx=\Int_{0}^{1}(g_4+\lambda g_3){\overline y}\, dx.
\hfill & \cr}\right.
\label{e25}
\end{equation}
Adding $(\ref{e25})_1$ and $(\ref{e25})_2$, we introduce a sesquilinear form ${\cal B}: [W_{a}^{1}(0,1)\times H_{0, a}^{1}(0,1)]^2\rightarrow
\C$ given by
$$
{\cal B}((u, v), (w, y))=\lambda^{2}\Int_{0}^{1}(u {\overline w}+v {\overline y})\, dx
+ \Int_{0}^{1}a(x)(u_x{\overline w_x}+v_x{\overline y_x})\, dx
+\alpha \Int_{0}^{1}(v{\overline w}+u{\overline y})\, dx +(\lambda {\tilde\zeta}+\beta) u(1){\overline w(1)},
$$
and a continuous antilinear functional ${\cal L}: W_{a}^{1}(0,1)\times H_{0, a}^{1}(0,1)\rightarrow \C$ where
$$
{\cal L}(w, y)= \Int_{0}^{1}(g_2+\lambda g_1){\overline w}\, dx-\zeta
\Int_{-\infty}^{+\infty} \Frac{\vartheta(\varsigma)}{\varsigma^{2}+\omega+\lambda}g_5(\varsigma)\,  d\varsigma{\overline w(1)}
+{\tilde\zeta} g_1(1) {\overline w(1)}+\Int_{0}^{1}(g_4+\lambda g_3){\overline y}\, dx,
$$
satisfying
\begin{equation}
{\cal B}((u, v), (w, y))={\cal L}(w, y), \label{e26}
\end{equation}
the sesquilinear form ${\cal B}(., .)$ is a bounded since for any
$(u, v), (w, y)\in W_{a}^{1}(0,1)\times H_{0, a}^{1}(0,1)$.
$$
\matrix{{\cal B}((u, v), (w, y))\leq \lambda^2 \|u\|_{L^2(0, 1)}\|w\|_{L^2(0, 1)}+\lambda^2 \|v\|_{L^2(0, 1)}\|y\|_{L^2(0, 1)}\hfill & \cr
+\|\sqrt{a(x)}u_x\|_{L^2(0, 1)}\|\sqrt{a(x)}w_x\|_{L^2(0, 1)}
+\|\sqrt{a(x)}v_x\|_{L^2(0, 1)}\|\sqrt{a(x)}y_x\|_{L^2(0, 1)}\hfill & \cr
+|\alpha|\|v\|_{L^2(0, 1)}\|w\|_{L^2(0, 1)}+|\alpha|\|u\|_{L^2(0, 1)}\|y\|_{L^2(0, 1)}
+\lambda {\tilde\zeta}|u(1)|^2\hfill & \cr
\leq M \|(u, v)\|_{W_{a}^{1}(0,1)\times H_{0, a}^{1}(0,1)}\|(w, y)\|_{W_{a}^{1}(0,1)\times H_{0, a}^{1}(0,1)},\hfill & \cr}
$$
and is coercive because $\forall (u, v)\in W_{a}^{1}(0,1)\times H_{0, a}^{1}(0,1)$
$$
\matrix{{\cal B}((u, v), (u, v))&\geq &\lambda^2 (\|u\|_{L^2(0, 1)}^2+\|v\|_{L^2(0, 1)}^2)+
\|\sqrt{a(x)}u_x\|_{L^2(0, 1)}^2+\|\sqrt{a(x)}v_x\|_{L^2(0, 1)}^2\hfill  \cr
&& +2\alpha\Re\Int_{0}^{1}u{\overline v}\, dx+\lambda {\tilde\zeta} |u(1)|^2\hfill  \cr
&\geq & c \|u\|_{W_{a}^{1}(0,1)}\|v)\|_{H_{0, a}^{1}(0,1)},\hfill  \cr}
$$
for $\alpha$ small enough.
Therefore, Lax-Milgram says that system (\ref{e26}) has a unique solution
$(u, v)\in W_{a}^{1}(0,1)\times H_{0, a}^{1}(0,1)$.

Now taking $(w, y)=(w, 0)$ with $w\in {\cal D}(0, 1)$ in (\ref{e26}), we obtain
\begin{equation}
\lambda^{2}u -(a(x)u_x)_x+\alpha v=g_2+\lambda g_1.
\label{mi1}
\end{equation}
Due to the fact that $u\in W_{a}^{1}(0,1)$ we get $(a(x)u_x)_x\in L^{2}(0, 1)$, and we deduce that
$u\in H_{a}^{2}(0,1)\cap W_{a}^{1}(0,1)$.

Similarly taking $(w, y)=(0, y)$ with $y\in {\cal D}(0, 1)$ in (\ref{e26}), we obtain
\begin{equation}
\lambda^{2}v -(a(x)v_x)_x+\alpha u=g_4+\lambda g_3,
\label{mi2}
\end{equation}
and we deduce that $v\in H_{a}^{2}(0,1)\cap H_{0, a}^{1}(0,1)$.

Multiplying both sides of the conjugate of equalities (\ref{mi1}) and (\ref{mi2}) by  $w\in W_{a}^{1}(0,1)$
and $y\in H_{0, a}^{1}(0,1)$, integrating by parts on $(0, 1)$, and comparing with  (\ref{e26}) we get
$$
\matrix{ (a(x)u_x)(1){\overline w}(1)
+(\varrho\lambda(\lambda+\omega)^{\tau-1} +\beta)u(1){\overline w}(1)\hfill & \cr
+\zeta\Int_{-\infty}^{+\infty} \Frac{\vartheta(\varsigma)}{\varsigma^{2}+\omega+\lambda}g_5(\varsigma)  \, d\varsigma{\overline w}(1)\
-\varrho(\lambda+\omega)^{\tau-1}  g_1(1){\overline w}(1)=0.\hfill & \cr}
$$
Consequently, defining ${\tilde u}= \lambda u- g_1$ and $\varphi$ by $(\ref{e18})_5$, we
deduce that
$$
\beta u(1)+(a(x)u_x)(1)+\zeta\Int_{-\infty}^{+\infty}\vartheta(\varsigma)\varphi(\varsigma)\, d\varsigma=0.
$$
In order to complete the existence of $U\in D({\cal A})$, we need to prove $\varphi$ and $|\varsigma|\varphi\in L^{2}(-\infty, \infty)$.
From $(\ref{e18})_5$, we get
$$
\Int_{\R}|\varphi(\varsigma)|^2\, d\varsigma\leq 3\Int_{\R}\Frac{|g_5(\varsigma)|^2}{(\varsigma^{2}+\omega+\lambda)^2}  \, d\varsigma
+3 (\lambda^2|u(1)|^2+|g_1(1)|^2)\Int_{\R}\Frac{|\varsigma|^{2\tau-1}}{(\varsigma^{2}+\omega+\lambda)^2}  \, d\varsigma.
$$
Using Proposition \ref{th2}, it easy to see that
$$
\Int_{\R}\Frac{|\varsigma|^{2\tau-1}}{(\varsigma^{2}+\omega+\lambda)^2}  \, d\varsigma=(1-\tau)\Frac{\pi}{\sin \tau\pi}(\lambda+\omega)^{\tau-2}.
$$
On the other hand, using the fact that $g_5\in L^2(\R)$, we obtain
$$
\Int_{\R}\Frac{|g_5(\varsigma)|^2}{(\varsigma^{2}+\omega+\lambda)^2}  \, d\varsigma\leq \Frac{1}{(\omega+\lambda)^2}
\Int_{\R}|g_5(\varsigma)|^2  \, d\varsigma< +\infty.
$$
It follows that $\varphi\in L^{2}(\R)$. Next, using $(\ref{e18})_5$, we get
$$
\Int_{\R}|\varsigma\varphi(\varsigma)|^2\, d\varsigma\leq 3\Int_{\R}\Frac{|\varsigma|^{2}|g_5(\varsigma)|^2}{(\varsigma^{2}+\omega+\lambda)^2}  \, d\varsigma
+3 (\lambda^2|u(1)|^2+|g_1(1)|^2)\Int_{\R}\Frac{|\varsigma|^{2\tau+1}}{(\varsigma^{2}+\omega+\lambda)^2}  \, d\varsigma.
$$
Using again Proposition \ref{th2}, it easy to see that
$$
\Int_{\R}\Frac{|\varsigma|^{2\tau+1}}{(\varsigma^{2}+\omega+\lambda)^2}  \, d\varsigma=\tau\Frac{\pi}{\sin \tau\pi}(\lambda+\omega)^{\tau-1}.
$$
Now, using the fact that $g_5\in L^2(\R)$, we obtain
$$
\Int_{\R}\Frac{|\varsigma|^{2}|g_5(\varsigma)|^2}{(\varsigma^{2}+\omega+\lambda)^2}  \, d\varsigma\leq \Frac{1}{(\omega+\lambda)}
\Int_{\R}|g_5(\varsigma)|^2  \, d\varsigma< +\infty.
$$
It follows that $|\varsigma|\varphi\in L^{2}(\R)$. Finally, since $\varphi\in L^{2}(\R)$, we get
$$
-(\varsigma^{2}+\omega)\varphi+{\tilde u}(1)\vartheta(\varsigma)=\lambda\varphi(\varsigma)-g_5(\varsigma)\in L^{2}(\R).
$$
Then $U\in D({\cal A})$ and
Therefore, the operator $\lambda I- {\cal A}$ is surjective for any $\lambda > 0$.

\hfill$\Box$\\


\subsection{Strong stability of the system}
In this part, we use a general criteria of Theorem \ref{thm14} to show
the strong stability of the $C_0$-semigroup $e^{t{\cal A}}$ associated to the wave system $(P)$
in the absence of the compactness of the resolvent of ${\cal A}$.

To state and prove our stability results, we need some results from semigroup theory.
\begin{theorem}[\cite{arba}]
Let ${\cal A}$ be the generator of a uniformly bounded
$C_0$-semigroup $\{S(t)\}_{t\geq 0}$ on a Hilbert space ${\cal X}$.
If:
\begin{itemize}
\item[(i)] ${\cal A}$ does not have eigenvalues on $i\R$.
\item[(ii)] The intersection of the spectrum $\sigma({\cal A})$ with $i\R$ is at most a countable set,
\end{itemize}
then the semigroup $\{S(t)\}_{t\geq 0}$ is asymptotically stable,
i.e, $\|S(t)z\|_{{\cal X}}\rightarrow 0$ as $t\rightarrow \infty$
for any $z\in {\cal X}$.
\label{thm14}
\end{theorem}
Our main result is the following theorem:
\begin{theorem}
The $C_0$-semigroup $e^{t{\cal A}}$ is strongly stable in ${\cal H}$; i.e, for all
$U_0\in {\cal H}$, the solution of (\ref{e14}) satisfies
$$
\Lim_{t\rightarrow \infty}\|e^{t{\cal A}}U_0\|_{\cal H}=0.
$$
\label{thmba}
\end{theorem}
For the proof of Theorem \ref{thmba}, we need the following two lemmas.
\begin{lemma}
${\cal A}$ does not have eigenvalues on $i\R$.
\label{lem10}
\end{lemma}
{\bf Proof.}\\
We will argue by contraction. Let $U\in D({\cal A})$ and let $\lambda\in \R$, such that
$$
{\cal A}U=i\lambda U.
$$
Then, we get

\begin{equation}
\left\{\matrix{i\lambda u-{\tilde u}=0,  \hfill & \cr
i\lambda {\tilde u}-(a(x)u_x)_x+\alpha v=0, \hfill & \cr
i\lambda v-{\tilde v}=0\hfill & \cr
i\lambda {\tilde v}-(a(x)v_x)_x+\alpha u=0, \hfill & \cr
i\lambda \varphi+(\varsigma^{2}+\omega)\varphi-{\tilde u}(1)\vartheta(\varsigma)=0.\hfill & \cr}\right.
\label{06z6}
\end{equation}
$\bullet${\bf Case 1}: If $\lambda\not=0$, then, from (\ref{e17})
we have
\begin{equation}
\varphi\equiv 0.
\label{e28}
\end{equation}
From $(\ref{06z6})_3$, we have
\begin{equation}
{\tilde u}(1)=0.
\label{e29}
\end{equation}
Hence, from $(\ref{06z6})_1$ we obtain
\begin{equation}
u(1)=0 \hbox{ and }  u_x(1)=0.
\label{e30}
\end{equation}
Eliminating ${\tilde u}$ and ${\tilde v}$ in equations $(\ref{06z6})_1$ and $(\ref{06z6})_3$ in equations $(\ref{06z6})_2$
and $(\ref{06z6})_4$, we obtain the following system
\begin{equation}
\left\{\matrix{
\lambda^2 u+(a(x)u_x)_x-\alpha v=0, \hfill & \cr
\lambda^2 v+(a(x)v_x)_x-\alpha u=0, \hfill & \cr
u(1)=u_{x}(1)=v(1)=0, \hfill & \cr
\left\{\matrix{u(0)=v(0)=0 \hfill &\hbox{ if } m_a\in [0, 1),\hfill \cr
(a(x) u_x)(0)=(a(x) v_x)(0)=0 \hfill &\hbox{ if } m_a\in [1, 2).\hfill \cr}\right.\hfill & \cr
}\right.
\label{v20}
\end{equation}
On the other hand, multiplying $(\ref{v20})_1$ by  ${\overline v}, (\ref{v20})_2$ by  ${\overline u}$ and
using the boundary condition $(\ref{v20})_3$, we get
\begin{equation}
\Int_{0}^{1}|u|^2\, dx=\Int_{0}^{1}|v|^2\, dx.
\label{v21}
\end{equation}
Multiplying equation $(\ref{v20})_1$ by $\overline{u}$, using Green formula, (\ref{e30})
and the boundary conditions, we get
\begin{equation}
\lambda^{2}\Int_{0}^{1}|u|^2\, dx-\Int_{0}^{1}a(x)|u_x|^2\, dx-\alpha\Int_{0}^{1}v\overline{u}\, dx=0.
\label{ee88}
\end{equation}
Multiplying equation $(\ref{v20})_1$ by $x\overline{u}_x$, we get
\begin{equation}
\lambda^{2}\Int_{0}^{1}x u \overline{u}_x\, dx+\Int_{0}^{1}x\overline{u}_x (a(x)u_x)_x\, dx
-\alpha\Int_{0}^{1}x v\overline{u}_x\, dx=0.
\label{ee89}
\end{equation}
$U\in D({\cal A})$, then the regularity is sufficiently for applying an integration on the second integral
in the left hand side in equation (\ref{ee89}). Then we obtain
\begin{equation}
\Frac{\lambda^{2}}{2}\Int_{0}^{1}x \Frac{d}{dx}|u|^2\, dx-\Int_{0}^{1}a(x)|u_x|^2\, dx
-\Frac{1}{2}\Int_{0}^{1}x a(x)\Frac{d}{dx}|u_x|^2\, dx-\alpha\Re\Int_{0}^{1}x v\overline{u_x}\, dx=0.
\label{ee90}
\end{equation}
Using Green formula, Proposition \ref{propo4}-(ii) and the boundary conditions, we get
\begin{equation}
\lambda^{2}\Int_{0}^{1}|u|^2\, dx+\Int_{0}^{1}(a(x)-xa'(x))|u_x|^2\, dx+2\alpha\Re\Int_{0}^{1}x v\overline{u_x}\, dx=0.
\label{ee91}
\end{equation}
Multiplying equations (\ref{ee88}) by $-m_a/2$, and tacking the sum of this equation and (\ref{ee91}), we
get
\begin{equation}
\matrix{\Frac{2-m_a}{2}\lambda^{2}\Int_{0}^{1}|u|^2\, dx+\Int_{0}^{1}\left(a(x)-xa'(x)+\Frac{m_a}{2}a(x)\right)|u_x|^2\, dx\hfill & \cr
+2\alpha\Re\Int_{0}^{1}x v\overline{u_x}\, dx+\alpha\Frac{m_a}{2}\Int_{0}^{1} v\overline{u}\, dx=0.\hfill & \cr}
\label{ee92}
\end{equation}
By definition of $m_a$, we have
$$
(2-m_a)a(x)\leq 2(a(x)-xa'(x))+m_a a(x).
$$
Then using the Cauchy-Schwartz and Poincar\'e's inequalities,
we deduce from (\ref{ee92}) and (\ref{v21}) that there exists a positive constant $C >0$,
$$
\Int_{0}^{1}a(x)|u_x|^2\, dx\leq \alpha C \Int_{0}^{1}a(x)|u_x|^2\, dx.
$$
which yields u=0 for $\alpha$ small enough. It then follows from (\ref{v21}) that $v=0$, and from $(\ref{06z6})_1$ and
$(\ref{06z6})_3$ that ${\tilde u}={\tilde v}= 0$.

Consequently, we obtain $U = 0$, which contradict
the hypothesis $U\not= 0$. The proof has been completed.

\noindent
$\bullet${\bf Case 2}: Otherwise, if $\lambda =0$, the system (\ref{06z6}) becomes
\begin{equation}
\left\{\matrix{{\tilde u}={\tilde v}=0,  \hfill & \cr
(a(x) u_x)_x-\alpha v=0, \hfill & \cr
(a(x) v_x)_x-\alpha u=0, \hfill & \cr
(\varsigma^{2}+\omega)\varphi-{\tilde u}(1)\vartheta(\varsigma)=0.\hfill & \cr}\right.
\label{e27qqq}
\end{equation}
From $(\ref{e27qqq})_1$ and $(\ref{e27qqq})_4$ , we have
\begin{equation}
\varphi \equiv 0.
\label{e27qqqr}
\end{equation}
Multiplying equation $(\ref{e27qqq})_2$ by $\overline{u}$, $(\ref{e27qqq})_3$ by $\overline{v}$,
using Green formula and the boundary conditions, we get
\begin{equation}
\Int_{0}^{1}a(x)[|u_x|^2+|v_x|^2]\, dx+\beta|u(1)|^2+\alpha\Int_{0}^{1}v{\overline u}\, dx+\alpha\Int_{0}^{1}u{\overline v}\, dx=0.
\label{elm}
\end{equation}
which yields $u_x = v_x = 0$ for $\alpha$ small enough.
Moreover, if $m_a\in [1, 2)$, then $u(1)=0$. Hence
$$
u=v=0.
$$
if $m_a\in [0, 1)$, then $u(0)=v(0)=0$.
Hence $u=v\equiv 0$.
and consequently, we obtain $U = 0$, which contradict
the hypothesis $U\not= 0$. The proof has been completed.

\begin{lemma}
We have
$$
\matrix{i\R\subset \rho({\cal A}) \hbox{ if }\omega\not=0, \hfill &\cr
i\R^*\subset \rho({\cal A}) \hbox{ if }\omega=0,\hfill &\cr}
$$
\label{lem11}
\end{lemma}
where $\R^*=\R-\{0\}$.\\
\noindent
{\bf Proof.}\\
\noindent
$\bullet${\bf Case 1}: $\lambda\not=0$.\\
We will prove that the operator $i\lambda I-{\cal A}$ is surjective for $\lambda\not=0$. For this purpose, let
$G=(g_1, g_2, g_3, g_4, g_5)^{T}\in {\cal H}$, we seek $X=(u, {\tilde u}, v, {\tilde v}, \varphi)^{T}\in D({\cal A})$ solution of the following equation
\begin{equation}
(i\lambda I-{\cal A})X=G.
\label{e35eee}
\end{equation}
Equivalently, we have
\begin{equation}
\left\{\matrix{i\lambda u-{\tilde u}=g_1,  \hfill & \cr
i\lambda {\tilde u}-(a(x)u_x)_x+\alpha v=g_2, \hfill & \cr
i\lambda v-{\tilde v}=g_3\hfill & \cr
i\lambda {\tilde v}-(a(x)v_x)_x+\alpha u=g_4, \hfill & \cr
i\lambda \varphi+(\varsigma^{2}+\omega)\varphi-{\tilde u}(1)\vartheta(\varsigma)=g_5.\hfill & \cr}\right.
\label{e35e}
\end{equation}
Inserting $(\ref{e35e})_1$ and $(\ref{e35e})_3$ into $(\ref{e35e})_2$ and $(\ref{e35e})_4$, we get
\begin{equation}
\left\{\matrix{-\lambda^2 u-(a(x)u_x)_x+\alpha v=(g_2+i\lambda g_1), \hfill & \cr
-\lambda^2 v-(a(x)v_x)_x+\alpha u=(g_4+i\lambda g_3), \hfill & \cr}\right.
\label{mk}
\end{equation}
Solving system (\ref{mk}) is equivalent to finding $u\in H_a^{2}\cap W_{a}^{1}(0,1)$ and $v\in H_a^{2}\cap H_{0, a}^{1}(0,1)$
such that
\begin{equation}
\left\{\matrix{\Int_{0}^{1}(-\lambda^{2}u {\overline w}- (a(x) u_x)_x {\overline w}+\alpha v{\overline w})\, dx=
\Int_{0}^{1}(g_2+i\lambda g_1){\overline w}\, dx,\hfill & \cr
\Int_{0}^{1}(-\lambda^{2}u {\overline y}- (a(x) u_x)_x {\overline y}+\alpha u{\overline y})\, dx=
\Int_{0}^{1}(g_4+i\lambda g_3){\overline y}\, dx,\hfill & \cr}\right.
\label{mk1}
\end{equation}
for all $w\in W_{a}^{1}(0,1)$ and $y\in H_{0, a}^{1}(0, 1)$. Then, we get
\begin{equation}
\left\{\matrix{
\Int_{0}^{1}(-\lambda^{2} u {\overline w}+ a(x) u_x  {\overline w_x}+\alpha v{\overline w})\, dx
+(i\lambda {\tilde\zeta}+\beta)  u(1)\ {\overline w(1)}\, \hfill &\cr
=\Int_{0}^{1}(g_2+i\lambda g_1){\overline w}\, dx
-\zeta \Int_{-\infty}^{+\infty} \Frac{\vartheta(\varsigma)}{\varsigma^{2}+\omega+i\lambda}f_5(\varsigma)  \, d\varsigma{\overline w}(1)
+{\tilde\zeta}  g_1(1) {\overline w(1)},\hfill & \cr
\Int_{0}^{1}(-\lambda^{2}v {\overline y}+ a(x) v_x {\overline y_x}+\alpha u{\overline y})\, dx=
\Int_{0}^{1}(g_4+i\lambda g_3){\overline y}\, dx.\hfill & \cr}\right.
\label{mk2}
\end{equation}
We can rewrite (\ref{mk2}) as
\begin{equation}
{\cal B}((u, v), (w, y))=l(w, y),\quad \forall (w, y)\in W_{a}^{1}\times H_{0, a}^{1}(0, 1),
\label{e2xx5}
\end{equation}
where
$$
{\cal B}(u, v), (w, y)={\cal B}_1(u, v), (w, y)+{\cal B}_2(u, v), (w, y)
$$
with
$$
\left\{\matrix{{\cal B}_1(u, v), (w, y)=\Int_{0}^{1}  (a(x)(u_x {\overline w_x}+v_x {\overline y_x})
+\alpha (v{\overline w}+ u{\overline y}))\, dx+(\beta u(1)+i\lambda {\tilde\zeta}) u(1)\ {\overline w(1)},\hfill \cr
{\cal B}_2(u, v), (w, y)=-\Int_{0}^{1}\lambda^{2}(u {\overline w}+v {\overline y})\, dx,\hfill \cr}\right.
\leqno{(*)}
$$
and
$$
\matrix{l(w, y)=\Int_{0}^{1}(g_2+i\lambda g_1){\overline w}\, dx
-\zeta\Int_{-\infty}^{+\infty} \Frac{\vartheta(\varsigma)}{\varsigma^{2}+\omega+i\lambda}g_3(\varsigma)  \, d\varsigma\ {\overline w}(1)\hfill &\cr
+\varrho(i\lambda+\omega)^{\tau-1}g_1(1){\overline w}(1)+\Int_{0}^{1}(g_4+i\lambda g_3){\overline y}\, dx.
&\cr}
$$
Let $(W_{a}^{1}\times H_{0, a}^{1}(0 , 1))'$ be the dual space of $W_{a}^{1}\times H_{0, a}^{1}(0 , 1)$.
Let us define the following operators
$$
\matrix{B:W_{a}^{1}\times H_{0, a}^{1}(0 , 1)\rightarrow (W_{a}^{1}\times H_{0, a}^{1}(0 , 1))'\hfill &\cr
\quad (u, v)\mapsto B (u, u) &\cr}
\matrix{B_i:W_{a}^{1}\times H_{0, a}^{1}(0 , 1)\rightarrow (W_{a}^{1}\times H_{0, a}^{1}(0 , 1))'\quad i\in \{1, 2\}\hfill &\cr
\quad (u, v)\mapsto B_i (u, v) &\cr}
\leqno{(**)}
$$
such that
$$
\matrix{(B (u, v)) (w, y)={\cal B}((u, v), (w, y)), \ \forall (w, y)\in W_{a}^{1}\times H_{0, a}^{1}(0 , 1),\hfill &\cr
(B_i u) w={\cal B}_i(u, w), \ \forall (w, y)\in W_{a}^{1}\times H_{0, a}^{1}(0 , 1), i\in \{1, 2\}.\hfill &\cr}
\leqno{(***)}
$$
We need to prove that the operator $B$ is an isomorphism. For this aim, we divide the proof into three steps:

\noindent
{\bf Step 1.}
In this step, we want to prove that the operator $B_1$ is an isomorphism. For this aim, it is easy to see that
${\cal B}_1$ is sesquilinear, continuous form on $W_{a}^{1}\times H_{0, a}^{1}(0 , 1)$.
Furthermore
$$
\matrix{\Re {\cal B}_1((u, v), (u, v))&=&\|\sqrt{a}u_x\|_{2}^2+\|\sqrt{a}v_x\|_{2}^2
+\alpha\Int_{0}^{1}(u{\overline v}+v{\overline u})\, dx+\beta|u(1)|^2\hfill \cr
&&+\varrho\lambda\Re \left(i(i\lambda+\omega)^{\tau-1}\right)|u(1)|^2\hfill \cr
&\geq & c(\|\sqrt{a}u_x\|_{2}^2+\|\sqrt{a}v_x\|_{2}^2+\beta|u(1)|^2),\hfill \cr}
$$
where we have used the fact that
$$
\varrho\lambda\Re \left(i(i\lambda+\omega)^{\tau-1}\right)=\zeta\lambda^2\Int_{-\infty}^{+\infty}\Frac{\vartheta(\varsigma)^2}{\lambda^2+(\omega+\varsigma^2)^2}\, d\varsigma> 0.
$$
Thus ${\cal B}_1$ is coercive. Then, from $(**)$ and Lax-Milgram theorem, the operator $B_1$ is
an isomorphism.

\noindent
{\bf Step 2.}
In this step, we want to prove that the operator $B_2$ is compact. For this aim, from $(*)$ and $(***)$, we have
$$
|{\cal B}_2((u, v), (w, y))|\leq c\|(u, v)\|_{L^{2}(0, 1)}\|(w, y)\|_{L^{2}(0, 1)},
$$
and consequently, using the compact embedding from $W_{a}^{1}\times H_{0, a}^{1}(0 , 1)$
to $L^{2}(0, 1)\times L^{2}(0, 1)$
we deduce that $B_2$ is a compact operator.
Therefore, from the above steps, we obtain that the operator $B =B_1+ B_2$ is a Fredholm operator of index zero. Now,
following Fredholm alternative, we still need to prove that the operator $B$ is injective to obtain that the operator $B$
is an isomorphism.

\noindent
{\bf Step 3.}
Let $(u, v)\in ker(B)$, then
\begin{equation}
{\cal B}(u, v), (w, y))=0\quad \forall (w, y)\in W_{a}^{1}\times H_{0, a}^{1}(0 , 1).
\label{eaz}
\end{equation}
In particular for $(w, y)=(u, v)$, it follows that
$$
\matrix{\lambda^{2}(\|u\|_{L^2(0, 1)}^2+\|v\|_{L^2(0, 1)}^2)
-i\varrho\lambda(i\lambda+\omega)^{\tau-1}|u(1)|^2-\beta|u(1)|^2=\hfill & \cr
\qquad\|\sqrt{a(x)}u_x\|_{L^2(0, 1)}^{2}+\|\sqrt{a(x)}v_x\|_{L^2(0, 1)}^{2}
+\alpha\Int_{0}^{1}(v{\overline u}+v{\overline u})\, dx. & \cr}
$$
Hence, we have
\begin{equation}
u(1)=0.
\label{e22xx1}
\end{equation}
From (\ref{eaz}), we obtain
\begin{equation}
(a(x)u_x)(1)=0
\label{e23xx1}
\end{equation}
and then
\begin{equation}
\left\{\matrix{
\lambda^2 u+(a(x)u_x)_x-\alpha v=0, \hfill & \cr
\lambda^2 v+(a(x)v_x)_x-\alpha u=0, \hfill & \cr
u(1)=u_{x}(1)=v(1)=0, \hfill & \cr
\left\{\matrix{u(0)=v(0)=0 \hfill &\hbox{ if } m_a\in [0, 1),\hfill \cr
(a(x) u_x)(0)=(a(x) v_x)(0)=0 \hfill &\hbox{ if } m_a\in [1, 2).\hfill \cr}\right.\hfill & \cr
}\right.
\label{aeq44b}
\end{equation}
Then, according to Lemma \ref{lem10}, we deduce that $(u, v)=(0, 0)$ and consequently $Ker(B)=\{0\}$.
Finally, from Step 3 and Fredholm alternative, we deduce that the operator $B$ is isomorphism. It is easy
to see that the operator $l$ is a antilinear and continuous form on $W_{a}^{1}\times H_{0, a}^{1}(0 , 1)$. Consequently,
(\ref{e2xx5}) admits a unique solution $(u, v)\in W_{a}^{1}\times H_{0, a}^{1}(0 , 1)$.
By using the classical elliptic regularity, we deduce that $U\in D({\cal A})$ is a unique solution
of (\ref{e35eee}). Hence $i\lambda- {\cal A}$ is surjective for all $\lambda\in \R^*$.\\

\noindent
{\bf Case 2: $\lambda=0$ and $\omega\not=0$.} Using Lax-Milgram Lemma, we obtain the result.

Taking account of Lemmas \ref{lem10}, \ref{lem11}
and from Theorem \ref{thm14} the $C_0$-semigroup $e^{t{\cal A}}$ is strongly stable in ${\cal H}$.

\hfill$\Box$\\

\subsection{Optimal condition for strong stability of the system in the case $a(x)=x^{\gamma}$}
\begin{theorem}
The $C_0$-semigroup $e^{t\mathcal{A}}$ is strongly stable in $\mathcal{H}$ if and only if the coefficient $\alpha$
satisfies
$$
\alpha\not=\Frac{1}{2}\left(\Frac{2-\gamma}{2}\right)^2 (j_{\nu_{\gamma}, k}^2-j_{\nu_{\gamma}, m}^2),\quad
k, m\in \N,
\leqno{(C)}
$$
where $\nu_{\gamma}=|1-\gamma|/(2-\gamma)$ and $j_{\nu, 1}< j_{\nu, 2}< \ldots <j_{\nu, k}<\ldots$ denote the sequence of positive zeros
of the Bessel function of first kind and of order $\nu$.
\label{thmba1}
\end{theorem}
\begin{equation}
\left\{\matrix{
\lambda^2 u+(x^{\gamma}u_x)_x-\alpha v=0, \hfill & \cr
\lambda^2 v+(x^{\gamma}v_x)_x-\alpha u=0, \hfill & \cr
u(1)=u_{x}(1)=v(1)=0, \hfill & \cr
\left\{\matrix{u(0)=v(0)=0 \hfill &\hbox{ if } m_a\in [0, 1),\hfill \cr
(x^{\gamma} u_x)(0)=(x^{\gamma} v_x)(0)=0 \hfill &\hbox{ if } m_a\in [1, 2).\hfill \cr}\right.\hfill & \cr
}\right.
\label{v30}
\end{equation}
We consider only the case $\gamma\in [0, 1[$. The case $\gamma\in [1, 2[$ is similar.
Then $\phi=u+v$ and $\psi=u-v$ satisfy
\begin{equation}
\left\{
\matrix{
(\lambda^2-\alpha)\phi+(x^{\gamma}\phi_{x})_x=0, \hfill & \cr
(\lambda^2+\alpha)\psi+(x^{\gamma}\psi_{x})_x=0. \hfill & \cr}
\right.
\label{e31qqq}
\end{equation}
The solution of the equation (\ref{e31qqq}) is given by
$$
\left\{
\matrix{
\phi(x)=c_1 \Phi_{+}(x) + c_2\Phi_{-}(x),\hfill & \cr
\psi(x)={\tilde c}_1 \Phi_{++}(x) + {\tilde c}_2\Phi_{--}(x),\hfill & \cr}
\right.
$$
togheter with the boundary conditions
$$
\phi(0)=\phi(1)=\psi(0)=\psi(1)=0,\quad u_{x}(1)=0,
$$
where
\begin{equation}
\left\{
\matrix{
\Phi_{+}(x)=x^{\frac{1-\gamma}{2}}J_{\nu_{\gamma}}\left(\frac{2}{2-\gamma}\sqrt{\lambda^2-\alpha} x^{\frac{2-\gamma}{2}}\right),
\quad  \Phi_{-}(x)=x^{\frac{1-\gamma}{2}}J_{-\nu_{\gamma}}\left(\frac{2}{2-\gamma}\sqrt{\lambda^2-\alpha} x^{\frac{2-\gamma}{2}}\right),\hfill &\cr
\Phi_{++}(x)=x^{\frac{1-\gamma}{2}}J_{\nu_{\gamma}}\left(\frac{2}{2-\gamma}\sqrt{\lambda^2+\alpha} x^{\frac{2-\gamma}{2}}\right),\hfill &\cr
\Phi_{--}(x)=x^{\frac{1-\gamma}{2}}J_{-\nu_{\gamma}}\left(\frac{2}{2-\gamma}\sqrt{\lambda^2+2\alpha} x^{\frac{2-\gamma}{2}}\right),
\hfill &\cr}
\right.
\label{mmrr}
\end{equation}
where
\begin{equation}
J_\nu(y)=\Sum_{m=0}^{\infty}\Frac{(-1)^m}{m!\Gamma(m+\nu+1)}\left(\Frac{y}{2}\right)^{2m+\nu}
=\Sum_{m=0}^{\infty}c_{\nu, m}^{+}y^{2m+\nu},
\label{bess1}
\end{equation}
\begin{equation}
J_{-\nu}(y)=\Sum_{m=0}^{\infty}\Frac{(-1)^m}{m!\Gamma(m-\nu+1)}\left(\Frac{y}{2}\right)^{2m-\nu}
=\Sum_{m=0}^{\infty}c_{\nu, m}^{-}y^{2m-\nu},
\label{bess2}
\end{equation}
where $J_{\nu}$ and $J_{-\nu}$ are Bessel functions of the first kind of order $\nu$ and $-\nu$.

As $\phi(0)=\psi(0)=0$, then $c_2={\tilde c}_2=0$. As $u(x)=\Frac{1}{2}(\phi(x)+\psi(x))$, we deduce that
$$
\phi_{x}(1)=-\psi_{x}(1).
$$
Then
$$
\matrix{c_1\{(1-\gamma)J_{\nu_{\gamma}}\left(\frac{2}{2-\gamma}\sqrt{\lambda^2-\alpha}\right)-\sqrt{\lambda^2-\alpha}
J_{\nu_{\gamma}+1}\left(\frac{2}{2-\gamma}\sqrt{\lambda^2-\alpha}\right)\}\hfill &\cr
=-{\tilde c}_1\{(1-\gamma)J_{\nu_{\gamma}}\left(\frac{2}{2-\gamma}\sqrt{\lambda^2+\alpha}\right)-\sqrt{\lambda^2+\alpha}
J_{\nu_{\gamma}+1}\left(\frac{2}{2-\gamma}\sqrt{\lambda^2+\alpha}\right)\},\hfill &\cr}
$$
Moreover $\phi(1)=\psi(1)=0$. Then
$$
\left\{
\matrix{
c_1 J_{\nu_{\gamma}}\left(\Frac{2}{2-\gamma}\sqrt{\lambda^2-\alpha}\right)=0,\quad
{\tilde c}_1 J_{\nu_{\gamma}}\left(\Frac{2}{2-\gamma}\sqrt{\lambda^2+\alpha}\right)=0,\hfill & \cr
c_1\{(1-\gamma)J_{\nu_{\gamma}}\left(\frac{2}{2-\gamma}\sqrt{\lambda^2-\alpha}\right)-\sqrt{\lambda^2-\alpha}
J_{\nu_{\gamma}+1}\left(\frac{2}{2-\gamma}\sqrt{\lambda^2-\alpha}\right)\}\hfill &\cr
=-{\tilde c}_1\{(1-\gamma)J_{\nu_{\gamma}}\left(\frac{2}{2-\gamma}\sqrt{\lambda^2+\alpha}\right)-\sqrt{\lambda^2+\alpha}
J_{\nu_{\gamma}+1}\left(\frac{2}{2-\gamma}\sqrt{\lambda^2+\alpha}\right)\},\hfill &\cr
}
\right.
$$
If Bessel are zero then
$$
\Frac{2}{2-\gamma}\sqrt{\lambda^2-\alpha}=j_{\nu_{\gamma}, k} \hbox{ and }\Frac{2}{2-\gamma}\sqrt{\lambda^2+\alpha}=j_{\nu_{\gamma}, m}
$$
for some integers $k$ and $m$. Hence, eigenvalues on $i\R$ exist iff
$$
\alpha=\Frac{1}{2}\left(\Frac{2-\gamma}{2}\right)^2 (j_{\nu_{\gamma}, k}^2-j_{\nu_{\gamma}, m}^2).
$$
Hence, if condition $(C)$ is satisfied we deduce that $c_1 = 0$ or ${\tilde c}_1 = 0$ and consequently $u = v = 0$.

Therefore $U=0$. Consequently, $\mathcal{A}$ does not have purely imaginary eigenvalues.

\subsection{Lack of exponential stability}
This section will be devoted to the study of the lack of exponential decay of solutions associated with the system
$(P')$.
\begin{proposition}
The $C_0$-semigroup of contractions $S(t) = e^{{\cal A} t}$ associated with (\ref{e14}) is not exponentially stable.
\label{pro1}
\end{proposition}
{\bf Proof.}
Let $\mu_n$ be an eigenvalue of ${\cal K}u=-(au_x)_x$ in $H_{0, a}^{1}(0, 1)$ corresponding to the normalized
eigenfunction $e_n$, and
$$
U_n=\Frac{1}{\sqrt{2}}\left(0, 0, \frac{e_n}{i\sqrt{\mu_n}}, e_n, 0\right)^T.
$$
Then a straightforward computation gives
$$
\|U_n\|_{{\cal H}}=1, \qquad \|(i\sqrt{\mu_n}-{\cal A})U_n\|_{{\cal H}}^2=\Frac{\alpha^2}{2\mu_n}\rightarrow 0.
$$
This shows that the resolvent of ${\cal A}$ is not uniformly bounded on the imaginary axis. Following
{\bf\cite{12}} and {\bf\cite{13}}, the system $(P')$ is not uniformly and exponentially stable in the energy
space ${\cal H}$.\\

\noindent
{\bf Precise spectral analysis in the case $a(x)=x^{\gamma}$.}\\
We aim to show that an
infinite number of eigenvalues of ${\cal A}$ approach the imaginary
axis which prevents the system $(P)$ from being exponentially
stable. Indeed we first compute the characteristic equation that
gives the eigenvalues of  ${\cal A}$.
Let $\lambda$ be an eigenvalue of ${\cal A}$ with associated
eigenvector $U=(u, v, \varphi)^{T}$. We consider only the case $\gamma\in [0, 1[$. The case $\gamma\in [1, 2[$ is similar.
Then ${\cal A}U = \lambda U$ is
equivalent to
\begin{equation}
\left\{\matrix{\lambda u-{\tilde u}=0,  \hfill & \cr
\lambda {\tilde u}-(x^{\gamma}u_x)_x+\alpha v=0, \hfill & \cr
\lambda v-{\tilde v}=0\hfill & \cr
\lambda {\tilde u}-(x^{\gamma}v_x)_x+\alpha u=0, \hfill & \cr
(\lambda+\varsigma^2+\omega)\varphi-{\tilde u}(1)\vartheta(\varsigma)=0,
}\right.
\label{066}
\end{equation}
with boundary conditions
\begin{equation}
\left\{\matrix{u(0)=v(0)=v(1)=0, \hfill & \cr
(\beta +\varrho \lambda(\lambda+\omega)^{\tau-1})u(1)+ u_{x}(1)=0.\hfill &\cr}\right.
\label{e556aa}
\end{equation}
Inserting $(\ref{066})_1$ into $(\ref{066})_2$ and $(\ref{066})_3$ into $(\ref{066})_4$, we get
\begin{equation}
\left\{\matrix{
\lambda^2 u-(x^{\gamma}u_x)_x+\alpha v=0, \hfill & \cr
\lambda^2 v-(x^{\gamma}v_x)_x+\alpha u=0, \hfill & \cr
u(0)=v(0)=v(1)=0, \hfill & \cr
(\beta +\varrho \lambda(\lambda+\omega)^{\tau-1})u(1)+ u_{x}(1)=0.
}\right.
\label{e331cc}
\end{equation}
Let us set
\begin{equation}
\left\{\matrix{ \phi=u+v, \hfill \cr
 {\psi}=u-v.
\hfill &\cr}\right.
\label{e771}
\end{equation}
Then, we obtain
 \begin{equation}
 \left\{\matrix{(\lambda^2+\alpha)\phi-(x^{\gamma}\phi_x)_x=0,\hfill &\cr
(\lambda^2-\alpha)\psi-(x^{\gamma}\psi_x)_x=0.
\hfill &\cr}\right.
\label{e775}
 \end{equation}
The solution of equations (\ref{e775}) is given by
\begin{equation}
\left\{\matrix{\phi(x)=c_{1}\Phi_{+}+c_{-}\Phi_{-},\hfill &\cr
\psi(x)={\tilde c}_{1}\Phi_{++}+{\tilde c}_{-}\Phi_{--},\hfill &\cr}\right.
\end{equation}
where $\Phi_{+}, \Phi_{-}, \Phi_{++}$ and $\Phi_{--}$ are defined by
$$
\left\{\matrix{\Phi_{+}(x):=x^{\frac{1-\gamma}{2}}J_{\nu_\gamma}\left(\frac{2}{2-\gamma}i{\tilde\lambda}
x^{\frac{2-\gamma}{2}}\right),\hfill &\cr
\Phi_{-}(x):=x^{\frac{1-\gamma}{2}}J_{-\nu_\gamma}\left(\frac{2}{2-\gamma}i{\tilde\lambda}
x^{\frac{2-\gamma}{2}}\right)\hfill &\cr}\right.
$$
and
$$
\left\{\matrix{\Phi_{++}(x):=x^{\frac{1-\gamma}{2}}J_{\nu_\gamma}\left(\frac{2}{2-\gamma}i{\tilde{\tilde\lambda}}
x^{\frac{2-\gamma}{2}}\right),\hfill &\cr
\Phi_{--}(x):=x^{\frac{1-\gamma}{2}}J_{-\nu_\gamma}\left(\frac{2}{2-\gamma}i{\tilde{\tilde\lambda}}
x^{\frac{2-\gamma}{2}}\right),\hfill &\cr}\right.
$$
where
$$
\nu_\gamma=\Frac{1-\gamma}{2-\gamma}.
$$
Then
$$
\left\{\matrix{u(x)=\Frac{1}{2}(c_{1}\Phi_{+}+c_{-}\Phi_{-}+{\tilde c}_{1}\Phi_{++}+{\tilde c}_{-}\Phi_{--}),
\hfill &\cr
v(x)=\Frac{1}{2}(c_{1}\Phi_{+}+c_{-}\Phi_{-}-{\tilde c}_{1}\Phi_{++}-{\tilde c}_{-}\Phi_{--}).
\hfill &\cr
}\right.
$$
Then, using the series expansion of $J_{\nu_{\alpha}}$ and
$J_{-\nu_{\alpha}}$, one obtains
$$
\matrix{\Phi_{+}(x)=\Sum_{m=0}^{\infty}{\tilde c}_{\nu_\gamma, m}^{+}x^{1-\gamma+(2-\gamma)m},\quad
\Phi_{-}(x)=\Sum_{m=0}^{\infty}{\tilde c}_{\nu_\gamma, m}^{-}x^{(2-\gamma)m}\hfill &\cr
\Phi_{++}(x)=\Sum_{m=0}^{\infty}{\tilde e}_{\nu_\gamma, m}^{+}x^{1-\gamma+(2-\gamma)m},\quad
\Phi_{--}(x)=\Sum_{m=0}^{\infty}{\tilde e}_{\nu_\gamma, m}^{-}x^{(2-\gamma)m}\hfill &\cr}
$$
with
$$
\left\{\matrix{{\tilde c}_{\nu_\gamma, m}^{+}=c_{\nu_\gamma, m}^{+}\left(\Frac{2}{2-\gamma}i{\tilde\lambda}\right)^{2m+\nu_\gamma},\quad
{\tilde c}_{\nu_\gamma, m}^{-}=c_{\nu_\gamma, m}^{-}\left(\Frac{2}{2-\gamma}i{\tilde\lambda}\right)^{2m-\nu_\gamma}\hfill & \cr
{\tilde e}_{\nu_\gamma, m}^{+}=c_{\nu_\gamma, m}^{+}\left(\Frac{2}{2-\gamma}i{\tilde{\tilde\lambda}}\right)^{2m+\nu_\gamma},\quad
{\tilde e}_{\nu_\gamma, m}^{-}=c_{\nu_\gamma, m}^{-}\left(\Frac{2}{2-\gamma}i{\tilde{\tilde\lambda}}\right)^{2m-\nu_\gamma}\hfill & \cr}\right.
$$
Next one easily verifies that $\Phi_{+}, {\tilde\Phi}_{+}\in H_{0, a}^1(0, 1)$: indeed,
$$
\matrix{{\Phi}_{+}(x)\sim_0 {\tilde c}_{\nu_\gamma, 0}^{+}x^{1-\gamma},\quad
x^{\gamma/2}{\Phi}'_{+}(x)\sim_0 (1-\gamma){\tilde c}_{\nu_\gamma, 0}^{+}x^{-\gamma/2},\hfill \cr
{\Phi}_{-}(x)\sim_0 {\tilde c}_{\nu_\gamma, 0}^{-},\quad
x^{\gamma/2}{\Phi}'_{-}(x)\sim_0 (2-\gamma){\tilde c}_{\nu_\gamma, 0}^{-}x^{1-\gamma/2},\hfill \cr}
$$
where we have used the following relation
\begin{equation}
x J'_{\nu}(x)=\nu J_{\nu}(x)-x J_{\nu+1}(x).
\label{e7112}
\end{equation}
Hence, given $c_{-}={\tilde c}_{-}=0,  u(x)=\frac{1}{2}(c_{1}\Phi_{+}(x)+{\tilde c}_{1}\Phi_{++}(x))\in H_{0, a}^1(0, 1)$
and $v(x)=\frac{1}{2}(c_{1}\Phi_{+}(x)-{\tilde c}_{1}\Phi_{++}(x))\in H_{0, a}^1(0, 1)$
with the boundary conditions
$$
\left\{\matrix{v(1)=0, \hfill & \cr
(\beta +\varrho \lambda(\lambda+\omega)^{\tau-1})u(1)+ u_{x}(1)=0.\hfill & \cr}\right.
$$
Then
\begin{equation}
M\pmatrix{c_{1}\cr {\tilde c}_{1}\cr}=\pmatrix{0\cr 0\cr},
\label{mat1}
\end{equation}
where
$$
M=\pmatrix{(\beta +\varrho \lambda(\lambda+\omega)^{\tau-1}){\Phi}_{+}(1)+{\Phi}'_{+}(1)&
(\beta +\varrho \lambda(\lambda+\omega)^{\tau-1}){\Phi}_{++}(1)+{\Phi}'_{++}(1)\cr
{\Phi}_{+}(1)& -{\Phi}_{++}(1)\cr}
$$
System (\ref{e331cc}) admits a non trivial solution if and only if $det(M)=0$. i.e.,
if and only if the eigenvalues of ${\cal A}$ are roots of the function $f$ defined by
\begin{equation}
\matrix{f(\lambda)&=&
2(\beta +(1-\gamma)+\varrho \lambda(\lambda+\omega)^{\tau-1})
J_{\nu_\gamma}\left(\frac{2}{2-\gamma}i{\tilde\lambda}\right)J_{\nu_\gamma}\left(\frac{2}{2-\gamma}i{\tilde{\tilde\lambda}}\right)\hfill \cr
&&-i{\tilde\lambda}J_{1+\nu_\gamma}\left(\frac{2}{2-\gamma}i{\tilde\lambda}\right)J_{\nu_\gamma}\left(\frac{2}{2-\gamma}i{\tilde{\tilde\lambda}}\right)
-i{\tilde{\tilde\lambda}}J_{1+\nu_\gamma}\left(\frac{2}{2-\gamma}i{\tilde{\tilde\lambda}}\right)J_{\nu_\gamma}\left(\frac{2}{2-\gamma}i{\tilde\lambda}\right). \cr}
\label{sszz}
\end{equation}
Our purpose is to prove, thanks to Rouch\'e's Theorem,
that there is a subsequence of eigenvalues for which their real part
tends to $0$.

In the sequel, since ${\cal A}$ is dissipative, we study the
asymptotic behavior of the large eigenvalues $\lambda$ of ${\cal A}$
in the strip $S=\{\lambda\in \C: -\alpha_0\leq \Re(\lambda)\leq 0\}$, for some $\alpha_0 > 0$ large enough and for such $\lambda$,
we remark that $\Phi_{+}, \Phi_{-}$
remain bounded.
\begin{lemma}
The large eigenvalues of the dissipative operator ${\cal A}$ are simple and can
be split into two families $(\lambda_{j, k})_{k\in \Z, |k|\geq N}, j = 1, 2, (N\in  \N$, chosen large enough).
Moreover, the following asymptotic expansions for the eigenvalues hold:

\noindent
$\bullet$ If $\tau=1$, then
$$
\lambda_{1, k}=\left\{\matrix{\Frac{2-\gamma}{2} \left[\ln\sqrt{\Frac{\varrho-1}{\varrho+1}}+i \left(k+\Frac{(1-2\nu_{\gamma})}{4}\right)\pi\right]+O\left(\Frac{1}{k}\right) \hfill & \hbox{ if }\varrho>1  \hfill\cr
\Frac{2-\gamma}{2} \left[\ln\sqrt{\Frac{1-\varrho}{\varrho+1}}+i\left(k+\Frac{3-2\nu_{\gamma}}{4}\right)\pi\right]+O\left(\Frac{1}{k}\right)\hfill & \hbox{ if }\varrho<1  \hfill\cr}\right\},\quad k\in\Z.
$$
$$
\lambda_{2, k}=\Frac{2-\gamma}{2}i\left(k\pi+\frac{(1-2\nu_{\gamma})}{4}\pi-\frac{a_1}{k\pi}+\frac{(1-2\nu_{\gamma})a_1}{4k^2\pi}\right)
-\left(\Frac{2}{2-\gamma}\right)^3\frac{\varrho\alpha^2}{4(k\pi)^{2}}
+O\left(\frac{1}{k^{3}}\right).
$$

\noindent
$\bullet$ If $0<\tau< 1$, then
$$
\lambda_{1, k}=\Frac{2-\gamma}{2}i\left(k+\frac{(3-2\nu_{\gamma})}{4}\right)\pi
+\frac{\beta_1}{k^{1-\tau}}+\frac{{\tilde\beta}_1}{k^{1-\tau}}
+o\left(\frac{1}{k^{1-\tau}}\right), k\geq N, {\tilde\beta}_1\in
i\R.
$$
where
$$
\beta_1=-\left(\Frac{2}{2-\gamma}\right)^{\tau}\Frac{\varrho}{\pi^{1-\tau}}\cos((1-\tau)\frac{\pi}{2}).
$$
$$
\lambda_{2, k}=\Frac{2-\gamma}{2}i\left(k\pi+\frac{(1-2\nu_{\gamma})}{4}\pi-\frac{a_1}{k\pi}+\frac{(1-2\nu_{\gamma})a_1}{4k^2\pi}\right)
+\frac{\beta_2}{k^{3-\tau}}+\frac{{\tilde\beta}_2}{k^{3-\tau}}
+o\left(\frac{1}{k^{3-\tau}}\right), k\geq N, {\tilde\beta}_2\in
i\R,
$$
where
$$
\beta_2=-\left(\Frac{2}{2-\gamma}\right)^{4-\tau}\Frac{\varrho\alpha^2}{4\pi^{3-\tau}}\cos(1-\tau)\frac{\pi}{2}.
$$
$$
\lambda_{j, k}=\overline{\lambda_{j, -k}} \hbox{ if } k\leq -N,
$$
Moreover for all $|k|\geq N$, the eigenvalues $\lambda_{j, k}$ are
simple.
\label{ll12}
\end{lemma}

\noindent
{\bf Proof.}

\noindent
{\bf Step 1.}
We will use the following classical development (see
{\bf\cite{lebed}} p. 122, (5.11.6)): for all $\delta>0$, the
following development holds when $|\arg z|< \pi-\delta$:
\begin{equation}
\matrix{J_{\nu}(z)=\left(\Frac{2}{\pi z}\right)^{1/2}\left[\cos\left(z-\nu\frac{\pi}{2}-\frac{\pi}{4}\right)
-\Frac{(\nu-\Frac{1}{2})(\nu+\Frac{1}{2})}{2}\Frac{\sin\left(z-\nu\frac{\pi}{2}-\frac{\pi}{4}\right)}{z}\right.\hfill & \cr
\left.-\Frac{(\nu-\Frac{1}{2})(\nu+\Frac{1}{2})(\nu-\Frac{3}{2})(\nu+\Frac{3}{2})}{8}\Frac{\cos\left(z-\nu\frac{\pi}{2}-\frac{\pi}{4}\right)}{z^2}
+O\left(\frac{1}{|z|^3}\right)\right]. & \cr}
\label{029}
\end{equation}
Moreover, for $\lambda$ large enough, we have
\begin{equation}
\matrix{{\tilde\lambda}=\sqrt{\lambda^2+\alpha}&=&\lambda+\Frac{\alpha}{2\lambda}-\Frac{\alpha^2}{8\lambda^3}+\Frac{O(\alpha^3)}{\lambda^5},\hfill \cr
{\tilde{\tilde\lambda}}=\sqrt{\lambda^2-\alpha}&=&\lambda-\Frac{\alpha}{2\lambda}-\Frac{\alpha^2}{8\lambda^3}+\Frac{O(\alpha^3)}{\lambda^5},\hfill \cr
{\tilde\lambda}+{\tilde{\tilde\lambda}}&=&2\lambda-\Frac{\alpha^2}{4\lambda^2}+\Frac{O(\alpha^2)}{\lambda^4},\hfill \cr
{\tilde\lambda}-{\tilde{\tilde\lambda}}&=&\Frac{\alpha}{\lambda}+\Frac{O(\alpha^2)}{\lambda^4},\hfill \cr
\lambda(\lambda+\omega)^{\tau-1}&=&\lambda^{\tau}+\Frac{(\tau-1)\omega}{\lambda^{1-\tau}}
+\Frac{(\tau-1)(\tau-2)\omega^2}{2\lambda^{2-\tau}}+O\left(\Frac{1}{\lambda^{3-\tau}}\right).\hfill \cr
}
\label{sq2}
\end{equation}
Then
$$
\matrix{-i{\tilde\lambda}J_{1+\nu_\gamma}\left(\frac{2}{2-\gamma}i{\tilde\lambda}\right)J_{\nu_\gamma}\left(\frac{2}{2-\gamma}i{\tilde{\tilde\lambda}}\right)
-i{\tilde{\tilde\lambda}}J_{1+\nu_\gamma}\left(\frac{2}{2-\gamma}i{\tilde{\tilde\lambda}}\right)J_{\nu_\gamma}\left(\frac{2}{2-\gamma}i{\tilde\lambda}\right)\hfill \cr
=-i\lambda\left[J_{1+\nu_\gamma}\left(\frac{2}{2-\gamma}i{\tilde\lambda}\right)J_{\nu_\gamma}\left(\frac{2}{2-\gamma}i{\tilde{\tilde\lambda}}\right)
+J_{1+\nu_\gamma}\left(\frac{2}{2-\gamma}i{\tilde{\tilde\lambda}}\right)J_{\nu_\gamma}\left(\frac{2}{2-\gamma}i{\tilde\lambda}\right)\right]
\hfill \cr
-i\Frac{\alpha}{\lambda}\left[J_{1+\nu_\gamma}\left(\frac{2}{2-\gamma}i{\tilde\lambda}\right)J_{\nu_\gamma}\left(\frac{2}{2-\gamma}i{\tilde{\tilde\lambda}}\right)
-J_{1+\nu_\gamma}\left(\frac{2}{2-\gamma}i{\tilde{\tilde\lambda}}\right)J_{\nu_\gamma}\left(\frac{2}{2-\gamma}i{\tilde\lambda}\right)\right]\hfill \cr
+O\left(\frac{1}{\lambda^4}\right).\cr}
$$
Now, we set
$$
\matrix{a=Z-\nu_{\gamma}\frac{\pi}{2}-\frac{\pi}{4}=\frac{2}{2-\gamma}i{\tilde\lambda}-\nu_{\gamma}\frac{\pi}{2}-\frac{\pi}{4},\hfill &\cr
b={\tilde Z}-\nu_{\gamma}\frac{\pi}{2}-\frac{\pi}{4}=\frac{2}{2-\gamma}i{\tilde{\tilde\lambda}}-\nu_{\gamma}\frac{\pi}{2}-\frac{\pi}{4}.\hfill &\cr}
$$
Thus, from (\ref{029}), we have
\begin{equation}
\left\{\matrix{J_{1+\nu_\gamma}(Z)=\sin a+\Frac{{\tilde a}_1(\nu_{\gamma})}{Z}\cos a-\sin a \Frac{{\tilde a}_2(\nu_{\gamma})}{Z^2}+O\left(\Frac{1}{Z^3}\right),\hfill \cr
\Phi_{++}(1)=J_{\nu_\gamma}({\tilde Z})=\cos b-\sin b\Frac{a_1(\nu_{\gamma})}{{\tilde Z}}-\cos b \Frac{a_2(\nu_{\gamma})}{{\tilde Z}^2}+O\left(\Frac{1}{Z^3}\right),\hfill \cr
J_{1+\nu_\gamma}({\tilde Z})=\sin b+\Frac{{\tilde a}_1(\nu_{\gamma})}{{\tilde Z}}\cos b-\sin b \Frac{{\tilde a}_2(\nu_{\gamma})}{{\tilde Z}^2}+O\left(\Frac{1}{Z^3}\right),\hfill \cr
\Phi_{+}(1)=J_{\nu_\gamma}(Z)=\cos a-\sin a\Frac{a_1(\nu_{\gamma})}{Z}-\cos a \Frac{a_2(\nu_{\gamma})}{Z^2}+O\left(\Frac{1}{Z^3}\right),\hfill \cr}\right.
\label{sq1}
\end{equation}
where
$$
\matrix{a_1(\nu_{\gamma})=-\Frac{\cos\nu_{\gamma} \pi}{\pi}\Frac{\Gamma(\frac{3}{2}+\nu_{\gamma})\Gamma(\frac{3}{2}-\nu_{\gamma})}{2},\
 {\tilde a}_1(\nu_{\gamma})=\Frac{\cos\nu_{\gamma} \pi}{\pi}\Frac{\Gamma(\frac{5}{2}+\nu_{\gamma})\Gamma(\frac{1}{2}-\nu_{\gamma})}{2}\hfill &\cr
a_2(\nu_{\gamma})=\Frac{\cos\nu_{\gamma} \pi}{\pi}\Frac{\Gamma(\frac{5}{2}+\nu_{\gamma})\Gamma(\frac{5}{2}-\nu_{\gamma})}{8},\
 {\tilde a}_2(\nu_{\gamma})=-\Frac{\cos\nu_{\gamma} \pi}{\pi}\Frac{\Gamma(\frac{7}{2}+\nu_{\gamma})\Gamma(\frac{3}{2}-\nu_{\gamma})}{8}.\hfill &\cr}
$$
Let us start with the case $\tau=1$.
Inserting (\ref{sq1}) and (\ref{sq2}) in (\ref{sszz}) we get
\begin{equation}
f(\lambda)=-i\lambda\left(\Frac{2}{\pi Z}\right)^{1/2}\left(\Frac{2}{\pi {\tilde Z}}\right)^{1/2}{\tilde f}(\lambda),
\label{030}
\end{equation}
where
$$
\matrix{{\tilde f}(\lambda)=-\cosh(*)+i\varrho(i\sinh(*)+1)\hfill &\cr
+\Frac{1}{\lambda}\left[-\frac{a_1}{{\tilde r}}(1-i\sinh(*))+\frac{{\tilde a}_1}{{\tilde r}}(1+i\sinh(*))
+i{\tilde\beta}(1+i\sinh(*))+\frac{2ia_1\varrho}{{\tilde r}}\cosh(*)\right]\hfill &\cr
+\Frac{1}{\lambda^2}\left[\frac{i\varrho r^2\alpha^2}{2}+\left(\frac{a_2+{\tilde a}_2}{{\tilde r}^2}
+\frac{a_1{\tilde a}_1}{{\tilde r}^2}
+\frac{2i{\tilde\beta}a_1}{{\tilde r}}\right)\cosh(*)\right.\hfill &\cr
\left.+2i\left(-\frac{\varrho a_2}{{\tilde r}^2}(1+i\sinh(*))
+\frac{\varrho a_1^2}{2{\tilde r}^2}(1-i\sinh(*))
\right)\right]+O(\Frac{1}{\lambda^3}),\hfill &\cr}
$$
where
$$
r=\Frac{2}{2-\gamma},\ {\tilde r}=i r,\ {\tilde\beta}=\beta+(1-\gamma)
$$
and
$$
(*)=2r\lambda+i\nu_{\gamma}\pi.
$$
Then
\begin{equation}
{\tilde f}(\lambda)=f_0(\lambda)+\Frac{f_1(\lambda)}{\lambda}
+\Frac{f_2(\lambda)}{\lambda^{2}}+O\left(\frac{1}{\lambda^3}\right),
\label{027}
\end{equation}
where
\begin{equation}
f_0(\lambda)=-\cosh(*)+i\varrho(i\sinh(*)+1),\hspace{3cm}
\label{032}
\end{equation}
\begin{equation}
\hspace{2cm}\ \matrix{f_1(\lambda)=-\frac{a_1}{{\tilde r}}(1-i\sinh(*))+\frac{{\tilde a}_1}{{\tilde r}}(1+i\sinh(*))
+i{\tilde\beta}(1+i\sinh(*))\hfill &\cr
 +\frac{2ia_1\varrho}{{\tilde r}}\cosh(*),&\cr}
\label{032s}
\end{equation}
\begin{equation}
\matrix{f_2(\lambda)=\frac{i\varrho r^2\alpha^2}{2}+\left(\frac{a_2+{\tilde a}_2}{{\tilde r}^2}
+\frac{a_1{\tilde a}_1}{{\tilde r}^2}
+\frac{2i{\tilde\beta}a_1}{{\tilde r}}\right)\cosh(*)
+2i\left(-\frac{\varrho a_2}{{\tilde r}^2}(1+i\sinh(*))\right.\hfill &\cr
\left.+\frac{\varrho a_1^2}{2{\tilde r}^2}(1-i\sinh(*))
\right).&\cr}
\label{032d}
\end{equation}
Note that $f_0, f_1$ and $f_2$ remain bounded in the strip $-\alpha_0\leq
\Re(\lambda)\leq 0$.

\noindent
{\bf Step 2.} We look at the roots of $f_0$. From
(\ref{032}), $f_0$ has has two families of roots that we denote $\lambda_{1, k}^0$ and $\lambda_{2, k}^0$.
$$
f_0(\lambda)=0 \Leftrightarrow -\cosh(*)+i\varrho(i\sinh(*)+1)=0,
$$
i.e
$$
-(\varrho+1)e^{4r\lambda}+2i\varrho e^{-\nu_{\gamma}\pi i}e^{2r\lambda}+(\varrho-1)e^{-2\nu_{\gamma}\pi i}=0.
$$
This yield
$$
\left\{\matrix{e^{2r\lambda}=\Frac{\varrho-1}{\varrho+1}e^{-\nu_{\gamma}\pi i} \hbox{ or }\hfill &\cr
e^{2r\lambda}=i e^{-\nu_{\gamma}\pi i}\hfill &\cr}\right.
$$
and directly implies that
$$
\matrix{\lambda_{1, k}^0&=&\left\{\matrix{\Frac{2-\gamma}{2} \left[\ln\sqrt{\Frac{\varrho-1}{\varrho+1}}+i \left(k+\Frac{(1-2\nu_{\gamma})}{4}\right)\pi\right] \hfill & \hbox{ if }\varrho>1  \hfill\cr
\Frac{2-\gamma}{2} \left[\ln\sqrt{\Frac{1-\varrho}{\varrho+1}}+i\left(k+\Frac{3-2\nu_{\gamma}}{4}\right)\pi\right]\hfill & \hbox{ if }\varrho<1  \hfill\cr}\right\},\quad k\in\Z,\hfill\cr
\lambda_{2, k}^0&=&\Frac{2-\gamma}{2}i\left(k\pi+\frac{(1-2\nu_{\gamma})}{4}\pi\right).\hfill\cr}
$$
Using Rouch\'e's Theorem, we deduce that ${\tilde f}$ admits an infinity of simple roots in $S$ denoted by $\lambda_{1, k}$ and $\lambda_{2, k}$ for
$|k|\geq k_0$, for $k_0$ large enough, such that
\begin{equation}
\lambda_{1, k}=\lambda_{1, k}^0+\varepsilon_{1, k},
\label{030eaa}
\end{equation}
\begin{equation}
\lambda_{2, k}=\lambda_{2, k}^0+\varepsilon_{2, k}.
\label{030e}
\end{equation}
{\bf Step 3. Asymptotic behavior of $\varepsilon_{1, k}$.}
We consider only the case $\varrho> 1$. The case  $\varrho< 1$ is similar.
Using (\ref{030eaa}), we get
\begin{equation}
\matrix{\sinh(*)=i(\cosh \ell+2r\varepsilon_{1, k}\sinh \ell +2 r^2\varepsilon_{1, k}^2\cosh \ell+o(\varepsilon_{1, k}^2)),\hfill \cr
\cosh(*)=i(\sinh \ell+2r\varepsilon_{1, k}\cosh \ell +2 r^2\varepsilon_{1, k}^2\sinh \ell+o(\varepsilon_{1, k}^2)),\hfill \cr}
\label{ddff}
\end{equation}
where $\ell=\ln\frac{\varrho-1}{\varrho+1}$.
Substituting (\ref{ddff}) into (\ref{032}), using that
${\tilde f}(\lambda_{1, k})=0$, we get
$$
\varepsilon_{1, k}=O\left(\Frac{1}{k}\right).
$$

\noindent
{\bf Step 4. Asymptotic behavior of $\varepsilon_{2, k}$.}
Using (\ref{030e}), we get
\begin{equation}
\matrix{\sinh(*)=i(1+2r^2\varepsilon_{2, k}^2+o(\varepsilon_{2, k}^2)),\hfill \cr
\cosh(*)=2i(r\varepsilon_{2, k}+\frac{2}{3}r^3\varepsilon_{2, k}^3+o(\varepsilon_{2, k}^3)).\hfill \cr}
\label{030e1}
\end{equation}
Substituting (\ref{030e1}) into (\ref{032}), using that
${\tilde f}(\lambda_{2, k})=0$, we get
\begin{equation}
{\tilde f}(\lambda_{2, k})=
-2ir\varepsilon_{2, k}-\Frac{2a_1}{{\tilde r}(\frac{ik\pi}{r})}+O(\varepsilon_{2, k}^2)+O\left(\Frac{1}{k^2}\right)+O\left(\Frac{\varepsilon_{2, k}}{k}\right)=0.
\label{hh1}
\end{equation}
The previous equation has one solution
\begin{equation}
\varepsilon_{2, k}=-i\frac{a_1}{r k\pi}+O(\varepsilon_{2, k}^2)+O\left(\Frac{1}{k^2}\right)+O\left(\Frac{\varepsilon_{2, k}}{k}\right).
\label{hh2}
\end{equation}

\noindent
We can write
\begin{equation}
\lambda_{2, k}=\Frac{2-\gamma}{2}i\left(k\pi+\frac{(1-2\nu_{\gamma})}{4}\pi\right)-i\frac{a_1}{r k\pi}+{\tilde\varepsilon}_{2, k},
\label{hh3}
\end{equation}
where ${\tilde\varepsilon}_{2, k}=o(1/k)$.
Substituting (\ref{hh3}) into (\ref{027}), we get
\begin{equation}
\matrix{{\tilde f}(\lambda_{2, k})=
-2ir{\tilde\varepsilon}_{2, k}+\Frac{2i\varrho a_1^2}{(k\pi)^2}-\Frac{i\varrho r^4 \alpha^2 }{2(k\pi)^2}-\Frac{(1-2\nu_{\gamma})a_1}{2k^2\pi}
-\Frac{4i\varrho a_1^2 }{(k\pi)^2}+\Frac{2i\varrho a_1^2 }{(k\pi)^2}+
+O({\tilde\varepsilon}_{2, k}^2)+O\left(\Frac{1}{k^3}\right)\hfill \cr
+O\left(\Frac{{\tilde\varepsilon}_{2, k}}{k}\right)=0.\cr}
\label{mml1}
\end{equation}
The previous equation gives
\begin{equation}
{\tilde\varepsilon}_{2, k}=-\Frac{\varrho r^3 \alpha^2 }{4(k\pi)^2}+\Frac{(1-2\nu_{\gamma})a_1}{4rk^2\pi}i
+O({\tilde\varepsilon}_{2, k}^2)+O\left(\Frac{1}{k^3}\right)+O\left(\Frac{{\tilde\varepsilon}_{2, k}}{k}\right).
\label{hh5}
\end{equation}
From (\ref{hh5}) we have in that case
$|k|^{2}\Re\lambda_{2, k}\sim \upsilon$ with
$$
\upsilon=-\-\Frac{\varrho r^3 \alpha^2 }{4\pi^2}.
$$

\noindent
{\bf Case $0<\tau< 1$}
$$
\matrix{{\tilde f}(\lambda)=-\cosh(*)+\Frac{i\varrho}{\lambda^{1-\tau}}(i\sinh(*)+1)\hfill &\cr
+\Frac{1}{\lambda}\left[-\frac{a_1}{{\tilde r}}(1-i\sinh(*))+\frac{{\tilde a}_1}{{\tilde r}}(1+i\sinh(*))
+i{\tilde\beta}(1+i\sinh(*))\right]\hfill &\cr
+\Frac{1}{\lambda^{2-\tau}}\left[\frac{2ia_1\varrho}{{\tilde r}}\cosh(*)
+i\varrho(1-\tau)\omega(1+i\sinh(*))\right]\hfill &\cr
+\Frac{1}{\lambda^2}\left(\frac{a_2+{\tilde a}_2}{{\tilde r}^2}
+\frac{a_1{\tilde a}_1}{{\tilde r}^2}
+\frac{2i{\tilde\beta}a_1}{{\tilde r}}\right)\cosh(*)\hfill &\cr
+\Frac{2i}{\lambda^{3-\tau}}\left[\frac{\varrho r^2\alpha^2}{4}-\frac{\varrho a_2}{{\tilde r}^2}(1+i\sinh(*))
+\frac{\varrho a_1^2}{2{\tilde r}^2}(1-i\sinh(*))
+\frac{\varrho(\tau-1)\omega a_1}{{\tilde r}}\cosh(*)\right.\hfill &\cr
\left.+\frac{\varrho (\tau-1)(\tau-2)\omega^2}{4}(1+i\sinh(*))\right]
+O\left(\Frac{1}{\lambda^3}\right).\hfill &\cr}
$$
Then
\begin{equation}
{\tilde f}(\lambda)=f_0(\lambda)+\Frac{f_1(\lambda)}{\lambda^{1-\tau}}+\Frac{f_2(\lambda)}{\lambda}
+\Frac{f_3(\lambda)}{\lambda^{2-\tau}}+\Frac{f_4(\lambda)}{\lambda^{2}}+\Frac{f_5(\lambda)}{\lambda^{3-\tau}}+O\left(\frac{1}{\lambda^3}\right),
\label{027h}
\end{equation}
where
\begin{equation}
f_0(\lambda)=-\cosh(*),
\label{032h}
\end{equation}
\begin{equation}
f_1(\lambda)=i\varrho(i\sinh(*)+1),
\label{032dh}
\end{equation}
\begin{equation}
f_2(\lambda)=-\frac{a_1}{{\tilde r}}(1-i\sinh(*))+\frac{{\tilde a}_1}{{\tilde r}}(1+i\sinh(*))
+i{\tilde\beta}(1+i\sinh(*)),
\label{032dh1}
\end{equation}
\begin{equation}
f_3(\lambda)=\frac{2ia_1\varrho}{{\tilde r}}\cosh(*),
 +i\varrho(1-\tau)\omega(1+i\sinh(*)),
\label{032sh}
\end{equation}
\begin{equation}
f_4(\lambda)=\left(\frac{a_2+{\tilde a}_2}{{\tilde r}^2}+\frac{a_1{\tilde a}_1}{{\tilde r}^2}
+\frac{2i{\tilde\beta}a_1}{{\tilde r}}\right)\cosh(*),
\label{032dh2}
\end{equation}
\begin{equation}
\matrix{f_5(\lambda)=\frac{i\varrho r^2\alpha^2}{2}
+2i\left[-\frac{\varrho a_2}{{\tilde r}^2}(1+i\sinh(*))+\frac{\varrho a_1^2}{2{\tilde r}^2}(1-i\sinh(*))\right.\hfill &\cr
\left.+\frac{\varrho(\tau-1)\omega a_1}{{\tilde r}}\cosh(*)
+\frac{\varrho (\tau-1)(\tau-2)\omega^2}{4}(1+i\sinh(*))\right]. &\cr}
\label{032dh3}
\end{equation}
Note that $f_0,f_1, f_2, f_3, f_4$ and $f_5$ remain bounded in the strip $-\alpha_0\leq
\Re(\lambda)\leq 0$.

\noindent
{\bf Step 2.} We look at the roots of $f_0$. From
(\ref{032h}), $f_0$ has two families of roots that we denote
$\lambda_{1, k}^0$ and $\lambda_{2, k}^0$.
$$
f_0(\lambda)=0 \Leftrightarrow -\cosh(*)=0.
$$
Then
$$
-e^{4r\lambda}-e^{-2\nu_{\gamma}\pi i}=0.
$$
Hence
$$
\left\{\matrix{e^{2r\lambda}=-e^{-\nu_{\gamma}\pi i} \hbox{ or }\hfill &\cr
e^{2r\lambda}=i e^{-\nu_{\gamma}\pi i}\hfill &\cr}\right.
$$
$$
\matrix{\lambda_{1, k}^0&=&
\Frac{2-\gamma}{2} i\left(k+\Frac{3-2\nu_{\gamma}}{4}\right)\pi,\quad k\in\Z,\hfill\cr
\lambda_{2, k}^0&=&\Frac{2-\gamma}{2}i\left(k\pi+\frac{(1-2\nu_{\gamma})}{4}\pi\right),\quad k\in\Z.\hfill\cr}
$$
Using Rouch\'e's Theorem, we deduce that ${\tilde f}$ admits an infinity of simple roots in $S$ denoted by $\lambda_{1, k}$ and $\lambda_{2, k}$ for
$|k|\geq k_0$, for $k_0$ large enough, such that
\begin{equation}
\lambda_{1, k}=\Frac{2-\gamma}{2}i\left(k\pi+\frac{(3-2\nu_{\gamma})}{4}\pi\right)+\varepsilon_{1, k},
\label{030eq}
\end{equation}
\begin{equation}
\lambda_{2, k}=\Frac{2-\gamma}{2}i\left(k\pi+\frac{(1-2\nu_{\gamma})}{4}\pi\right)+\varepsilon_{2, k}.
\label{030eg}
\end{equation}
{\bf Step 3. Asymptotic behavior of $\varepsilon_{1, k}$.}
Using (\ref{030eq}), we get
\begin{equation}
\matrix{\sinh(*)=-i(1+2r^2\varepsilon_{1, k}^2+o(\varepsilon_{1, k}^2))\hfill \cr
\cosh(*)=-2i(r\varepsilon_{1, k}+\frac{2}{3}r^3\varepsilon_{1, k}^3+o(\varepsilon_{1, k}^3)).\hfill \cr}
\label{030e1q}
\end{equation}
Substituting (\ref{030e1q}) into (\ref{027h}), using that
${\tilde f}(\lambda_{1, k})=0$, we get
\begin{equation}
{\tilde f}(\lambda_k)=
2ir\varepsilon_{1, k}+\Frac{2i\varrho}{(\frac{ik\pi}{r})^{1-\tau}}+O\left(\Frac{1}{k}\right)+O\left(\Frac{\varepsilon_{1, k}^2}{k^{1-\tau}}\right)=0.
\label{mml4}
\end{equation}
The previous equation has one solution
\begin{equation}
\varepsilon_{1, k}=-\frac{\varrho}{r^{\tau} (k\pi)^{1-\tau}}(\cos(1-\tau)\frac{\pi}{2}-\sin(1-\tau)\frac{\pi}{2})
+O\left(\Frac{1}{k}\right)+O\left(\Frac{\varepsilon_{1, k}^2}{k^{1-\tau}}\right).
\label{055zaq}
\end{equation}
From (\ref{055zaq}) we have in that case
$|k|^{1-\tau}\Re\lambda_{1, k}\sim \beta_1$ with
$$
\beta_1=-\Frac{\varrho r^{-\tau}  }{\pi^{1-\tau}}\cos(1-\tau)\frac{\pi}{2}.
$$

{\bf Step 4. Asymptotic behavior of $\varepsilon_{2, k}$.}

Using (\ref{030eg}), we get
\begin{equation}
\matrix{\sinh(*)=i(1+2r^2\varepsilon_{2, k}^2+o(\varepsilon_{2, k}^2)),\hfill \cr
\cosh(*)=2i(r\varepsilon_{2, k}+\frac{2}{3}r^3\varepsilon_{2, k}^3+o(\varepsilon_{2, k}^3)).\hfill \cr}
\label{030e1g}
\end{equation}
Substituting (\ref{030e1g}) into (\ref{027h}), using that
${\tilde f}(\lambda_{2, k})=0$, we get
\begin{equation}
{\tilde f}(\lambda_{2, k})=
-2ir\varepsilon_{2, k}-\Frac{2a_1}{{\tilde r}(\frac{ik\pi}{r})}+O(\varepsilon_{2, k}^3)+O\left(\Frac{1}{k^2}\right)+O\left(\Frac{\varepsilon_{2, k}^2}{k^{1-\tau}}\right)=0.
\label{mml1x}
\end{equation}
The previous equation has one solution
\begin{equation}
\varepsilon_{2, k}=-i\frac{a_1}{r k\pi}+O(\varepsilon_{2, k}^3)+O\left(\Frac{1}{k^2}\right)+O\left(\Frac{\varepsilon_{2, k}^2}{k^{1-\tau}}\right).
\label{055zad}
\end{equation}
We can write
\begin{equation}
\lambda_{2, k}=\Frac{2-\gamma}{2}i\left(k\pi+\frac{(1-2\nu_{\gamma})}{4}\pi\right)-i\frac{a_1}{r k\pi}+{\tilde\varepsilon}_{2, k},
\label{rr22n}
\end{equation}
where ${\tilde\varepsilon}_{2, k}=o(1/k)$.
Substituting (\ref{rr22n}) into (\ref{027h}), we get
\begin{equation}
\matrix{{\tilde f}(\lambda_{2, k})=
-2ir\varepsilon_{2, k}-\Frac{i\varrho r^{5-\tau} \alpha^2 }{2i^{1-\tau}(k\pi)^{3-\tau}}-\Frac{(1-2\nu_{\gamma})a_1}{2k^2\pi}
+O\left(\Frac{\varepsilon_{2, k}^2}{k^{1-\tau}}\right)+O\left(\Frac{\varepsilon_{2, k}}{k^{2-\tau}}\right)\hfill \cr
+O(\varepsilon_{2, k}^2)+O\left(\Frac{1}{k^3}\right)=0.\cr}
\label{mml3}
\end{equation}
The previous equation gives
\begin{equation}
\varepsilon_{2, k}=-\Frac{\varrho r^{4-\tau} \alpha^2 }{4i^{1-\tau}(k\pi)^{3-\tau}}+\Frac{(1-2\nu_{\gamma})a_1}{4rk^2\pi}i
+O\left(\Frac{\varepsilon_{2, k}^2}{k^{1-\tau}}\right)+O\left(\Frac{\varepsilon_{2, k}}{k^{2-\tau}}\right)
+O(\varepsilon_{2, k}^2)+O\left(\Frac{1}{k^3}\right).
\label{rr23c}
\end{equation}
From (\ref{rr23c}) we have in that case
$|k|^{3-\tau}\Re\lambda_{2, k}\sim \beta_2$ with
$$
\beta_2=-\Frac{\varrho r^{4-\tau} \alpha^2 }{4\pi^{3-\tau}}\cos(1-\tau)\frac{\pi}{2}.
$$
Now, setting ${\tilde U}_k=(\lambda_{j, k}^0-{\cal A}) U_k$, where $U_k$ is
a normalized eigenfunction associated to $\lambda_{j, k}$. We then have
$$
\matrix{\|(\lambda_{2, k}^0-{\cal A})^{-1}\|_{{\cal L}({\cal H})}=\Sup_{U\in {\cal H}, U\not=0} \Frac{\|(\lambda_{2, k}^0-{\cal A})^{-1}U\|_{{\cal H}}}{\|U\|_{{\cal H}}}
&\geq &\Frac{\|(\lambda_{2, k}^0-{\cal A})^{-1}{\tilde U}_k\|_{{\cal H}}}{\|{\tilde U}_k\|_{{\cal H}}}\hfill \cr
&\geq &\Frac{\|U_k\|_{{\cal H}}}{\|(\lambda_{2, k}^0-{\cal A}) U_k\|_{{\cal H}}}.\hfill \cr}
$$
Hence, by Lemma \ref{ll12}, we deduce that
$$
\|(\lambda_{2, k}^0-{\cal A})^{-1}\|_{{\cal L}({\cal H})}\geq c \left\{\matrix{|k|^{2}\hfill & \hbox{ if }\tau=1,\hfill \cr
|k|^{3-\tau}\hfill & \hbox{ if }0<\tau <1.\hfill \cr}\right.
$$
So that, the semigroup $e^{t {\cal A}}$
is not exponentially stable. Thus the proof is complete.

\hfill$\Box$\\

\section{Polynomial Stability (for $\omega\not=0$)}
To prove polynomial decay, we use the following theorem.
\begin{theorem}[\cite{boto}]
Assume that ${\cal A}$ is the
generator of a strongly continuous semigroup of contractions $(e^{t{\cal A}})_{t\geq 0}$
on a Hilbert space ${\cal X}$. If $i\R\subset\varrho({\cal A})$. Then for a fixed $l> 0$ the following conditions are equivalent
\begin{itemize}
\item[1)] $\Sup_{\beta\in \R}\|(i\beta I-{\cal A})^{-1}\|_{{\cal L}({\cal X})}=O(|\beta|^{l})$.
\item[2)] $\|e^{t{\cal A}}U_0\|_{\cal X}\leq \Frac{C}{t^{\frac{1}{l}}} \|U_0\|_{D({\cal A})}\quad\forall t>0,\ U_0\in D({\cal A})$,
for some $C> 0$.
\end{itemize}
\label{thm2}
\end{theorem}

\begin{theorem}
The semigroup ${S_{\cal A}(t)}_{t\geq 0}$ associated with system $(P')$ is polynomially stable, i.e., there exists
a constant $C> 0$ such that
$$
{\cal E}(t)=\|S_{\cal A}(t)U_0\|_{\cal H}^2\leq \Frac{C}{ t^{\frac{2}{(3-\tau)} }  }\|U_0\|_{D({\cal A})}^2.
$$
\label{thm225}
\end{theorem}
{\bf Proof}\\
In section 3, we have proved that the first condition in Theorem \ref{thm2} is satisfied. Now,we need to show that
\begin{equation}
\Sup_{|\lambda|\geq 1}\Frac{1}{\lambda^{l}}\|(i\lambda I-{\cal A})^{-1}\|_{{\cal H}}< \infty,
\label{v35}
\end{equation}
where $l=3-\tau$.
We establish (\ref{v35}) by contradiction. So, if (\ref{v35}) is false, then there exist sequences
$(\lambda_n)_n\subset \R$ and $U_n= (u_n, {\tilde u}_n, v_n, {\tilde v}_n, \varphi_n)\in  D({\cal A})$
satisfying
\begin{equation}
\|U_n\|_{{\cal H}}=1\quad \forall n\geq 0,
\label{v37}
\end{equation}
\begin{equation}
\Lim_{n\rightarrow\infty}|\lambda_n|=\infty
\label{v38}
\end{equation}
and
\begin{equation}
\Lim_{n\rightarrow\infty}\lambda_n^{l}\|(i\lambda_n I-{\cal A})U_n\|\rightarrow 0,
\label{v39}
\end{equation}
which implies that
\begin{equation}
\left\{\matrix{\lambda^{l}(i\lambda u-{\tilde u})=g_1\rightarrow 0 \hbox{ in } W_{a}^1(0, 1),  \hfill & \cr
\lambda^{l}(i\lambda{\tilde u}-(a(x)u_x)_x+\alpha v)=g_2\rightarrow 0 \hbox{ in } L^2(0, 1), \hfill & \cr
\lambda^{l}(i\lambda v-{\tilde v})=g_3\rightarrow 0 \hbox{ in } H_{a, 0}^1(0, 1)\hfill & \cr
\lambda^{l}(i\lambda {\tilde v}-(a(x)v_x)_x+\alpha u)=g_4\rightarrow 0 \hbox{ in } L^2(0, 1), \hfill & \cr
\lambda^{l}(i\lambda \varphi+(\varsigma^{2}+\omega)\varphi-{\tilde u}(1)\vartheta(\varsigma))=g_5\rightarrow 0 \hbox{ in } L^{2}(-\infty, +\infty).\hfill & \cr}\right.
\label{v4}
\end{equation}
For simplification, we denote $\lambda_n$ by $\lambda, U_n= (u_n, {\tilde u}_n, v_n, {\tilde v}_n, \varphi_n)$ by
$U= (u, {\tilde u}, v, {\tilde v}, \varphi)$ and
$H_n=(g_{1n}, g_{2n}, g_{3n}, g_{4n}, g_{5n})=\lambda_n^{l}(i\lambda_n I-{\cal A})U_n$ by
$G_n=(g_{1}, g_{2}, g_{3}, g_{4}, g_{5})$.
We will prove that\\ $\|U\|_{\cal H}=o(1)$ as a contradiction with (\ref{v37}).
Our proof is divided into several steps.

\noindent
$\bullet${\bf Step 1}
Taking the inner product of $\lambda^{l}(i\lambda I-{\cal A})U$ with $U$, we get
\begin{equation}
i\lambda \|U\|_{{\cal H}}^{2}-({\cal A}U, U)_{{\cal H}}=\Frac{o(1)}{\lambda^{l}}.
\label{v2}
\end{equation}
Using (\ref{e17}), we get
\begin{equation}
\zeta\int_{-\infty}^{+\infty}(\varsigma^{2}+\omega)|\varphi(\varsigma)|^{2}\, d\varsigma=-\Re({\cal A}U, U)=\Frac{o(1)}{\lambda^{l}}.
\label{v3}
\end{equation}
Now, from $(\ref{v4})_{5}$, we obtain
\begin{equation}
{\tilde u}(1)\vartheta(\varsigma)=(i\lambda+\varsigma^{2}+\omega)\varphi-\Frac{g_5(\varsigma)}{\lambda^l}.
\label{e40kk}
\end{equation}
By multiplying $(\ref{e40kk})$ by $(i\lambda+\varsigma^{2}+\omega)^{-2}|\varsigma|$, we get
\begin{equation}
(i\lambda+\varsigma^{2}+\omega)^{-2}{\tilde u}(1)\vartheta(\varsigma)|\varsigma|=(i\lambda+\varsigma^{2}+\omega)^{-1}|\varsigma|\varphi-(i\lambda+\varsigma^{2}+\omega)^{-2}|\varsigma|\Frac{g_5(\varsigma)}{\lambda^l}.
\label{e38kk}
\end{equation}
Hence, by taking absolute values of both sides of (\ref{e38kk}), integrating over the interval $]-\infty, +\infty[$ with
respect to the variable $\varsigma$ and applying Cauchy-Schwartz inequality, we obtain
\begin{equation}
{\cal R}|{\tilde u}(1)|\leq  \sqrt{2}{\cal P} \left(\Int_{-\infty}^{+\infty}\varsigma^{2}|\varphi|^{2}\, d\varsigma\right)^{\frac{1}{2}}
+ 2 \Frac{{\cal Q}}{\lambda^l}\left(\Int_{-\infty}^{+\infty}|g_5(\varsigma)|^{2}\, d\varsigma\right)^{\frac{1}{2}},
\label{e39kk}
\end{equation}
where
$$
{\cal R}=\left|\Int_{-\infty}^{+\infty}(i\lambda+\varsigma^{2}+\omega)^{-2}|\varsigma|\vartheta(\varsigma)\, d\varsigma\right|
=\frac{|1-2\tau|}{4}\Frac{\pi}{|\sin\frac{(2\tau+3)}{4}\pi|}|i\lambda+\omega|^{\frac{(2\tau-5)}{4}},
$$
$$
{\cal P}=\left(\Int_{-\infty}^{+\infty}(|\lambda|+\varsigma^{2}+\omega)^{-2}\, d\varsigma\right)^{\frac{1}{2}}=(\frac{\pi}{2})^{1/2}||\lambda|+\omega|^{-\frac{3}{4}},
$$
$$
{\cal Q}=\left(\Int_{-\infty}^{+\infty}(|\lambda|+\varsigma^{2}+\omega)^{-4}|\varsigma|^{2}\, d\varsigma\right)^{\frac{1}{2}}
=\left(\Frac{\pi}{16}||\lambda|+\omega|^{-\frac{5}{2}}\right)^{1/2}.
$$
Thus, by using the inequality $2PQ \leq P^2 + Q^2, P \geq 0, Q \geq 0$, again, we get
\begin{equation}
{\cal R}^2|{\tilde u}(1)|^2\leq  2 {\cal P}^2 \left(\Int_{-\infty}^{+\infty}(\varsigma^{2}+\omega)|\varphi|^{2}\, d\varsigma\right)
+ 4\Frac{{\cal Q}^2}{\lambda^{2l}}\left(\Int_{-\infty}^{+\infty}|g_5(\varsigma)|^{2}\, d\varsigma\right).
\label{e339kk}
\end{equation}
We deduce that
\begin{equation}
|{\tilde u}(1)|^2=\Frac{o(1)}{\lambda^{l-(1-\tau)}}+\Frac{o(1)}{\lambda^{2l+\tau}}.
\label{e41nkk}
\end{equation}
Then
\begin{equation}
|{\tilde u}(1)|=\Frac{o(1)}{\lambda^{\frac{l-(1-\tau)}{2}}}.
\label{mn20}
\end{equation}
So, from $(\ref{v4})_1$, we find
\begin{equation}
|u(1)|=\left|\Frac{v(1)}{i\lambda}+\Frac{g_1(1)}{i\lambda^{l+1}}\right|=\Frac{o(1)}{\lambda^{\frac{l-(1-\tau)}{2}+1}}.
\label{v5}
\end{equation}
Since $U\in D({\cal A})$ and using the boundary conditions $(\ref{e15})_3$, (\ref{v5}) and (\ref{v3}), we obtain
\begin{equation}
|a(1)u_x(1)|=\Frac{o(1)}{\lambda^{\frac{l-(1-\tau)}{2}}}.
\label{v33}
\end{equation}

\noindent
$\bullet${\bf Step 2} Now we use the classical multiplier method.
Let us introduce the following notation
$$
{\cal I}_{v}(x)= |\sqrt{a(x)}v_x(x)|^{2} +|{\tilde v}(x)|^{2}.
$$
For simplification, we set ${\tilde g}_1=\frac{g_1}{\lambda^l},\ {\tilde g}_2=\frac{g_2}{\lambda^l},\
{\tilde g}_3=\frac{g_3}{\lambda^l},\ {\tilde g}_4=\frac{g_4}{\lambda^l},\ {\tilde g}_5=\frac{g_5}{\lambda^l}$.

\begin{lemma}
We have that
\begin{equation}
\matrix{\Int_{0}^{1} \left[\left((a(x)-xa'(x))+\frac{m_a}{2}a(x)\right)|v_{x}|^{2}+\left(1-\frac{m_a}{2}\right)|{\tilde v}(x)|^{2}\right]\, dx\hfill &\cr
+2\alpha\Re\Int_{0}^{1}xu {\overline v_x}\, dx+\alpha\Frac{m_a}{2}\Re\Int_{0}^{1}u {\overline v}\, dx\hfill &\cr
=[x{\cal I}_{v}]_{0}^{1}+\frac{m_a}{2}[a(x)v_x\overline{v}]_{0}^{1}+R,\hfill &\cr}
\label{ee22ww}
\end{equation}
where
$$
R=2 \Re\Int_{0}^{1} x {\tilde g}_4  \overline{v}_x\, dx+2 \Re\Int_{0}^{1}x {\tilde v}\overline{{\tilde g}}_{3x}\, dx
+\frac{m_a}{2}\Int_{0}^{1} {\tilde v}\overline{{\tilde g}_{3}}\, dx+\frac{m_a}{2}\Int_{0}^{1}{\tilde g}_{4}\overline{v}dx.
$$
\label{l11}
\end{lemma}
{\bf Proof.}\\
To get (\ref{ee22ww}), let us multiply the equation $(\ref{v4})_4$ by $x \overline{v}_x$ Integrating on $(0, 1)$ we obtain
$$
i\lambda\Int_{0}^{1} {\tilde v} x \overline{v}_x\, dx-\Int_{0}^{1} (a(x)v_x)_x x \overline{v}_x\, dx
+\alpha\Int_{0}^{1}xu {\overline v_x}\, dx=\Int_{0}^{L}{\tilde g}_4 x \overline{v}_x\, dx
$$
or
$$
-\Int_{0}^{1} {\tilde v} x (\overline{i\lambda v_x})\, dx-\Int_{0}^{1} x (a(x)v_x)_x  \overline{v}_x\, dx
+\alpha\Int_{0}^{1}xu {\overline v_x}\, dx= \Int_{0}^{1}{\tilde g}_4 x \overline{v}_x\, dx.
$$
Since $i\lambda v_x={\tilde v}_x+{\tilde g}_{3x}$ taking the real part in the above equality, we get
$$
\matrix{-\Frac{1}{2}\Int_{0}^{1} x \Frac{d}{dx}|{\tilde v}|^{2}\, dx+\Frac{1}{2}\Int_{0}^{1}x a(x)\Frac{d}{dx}|v_{x}|^{2}\, dx
-[x a(x) |v_{x}|^{2} ]_{0}^{1}+\Int_{0}^{1}  a(x)|v_{x}|^{2}\, dx+\alpha\Re\Int_{0}^{1}xu {\overline v_x}\, dx
\hfill & \cr
=   \Re\Int_{0}^{1 }{\tilde v} x\overline{{\tilde g}}_{3x}\,dx +  \Re\Int_{0}^{1}{\tilde g}_4 x \overline{v}_x\, dx .\hfill & \cr}
$$
Performing an integration by parts we obtain
\begin{equation}\label{e35}
\matrix{\Int_{0}^{1}  [|\sqrt{a(x)} v_{x}|^{2}+|{\tilde v}(x)|^{2}]\, dx
-\Int_{0}^{1}  x a'(x)|v_{x}(x)|^{2}\, dx+2\alpha\Re\Int_{0}^{1}xu {\overline v_x}\, dx\hfill &\cr
\qquad\qquad=[x(|\sqrt{a(x)} v_{x}|^{2}+|{\tilde v}(x)|^{2})]_{0}^{1}+R_{1}, &\cr}
\end{equation}
where
$$
R_{1}=2 \Re\Int_{0}^{1} x {\tilde g}_4  \overline{v}_x\, dx+2 \Re\Int_{0}^{1}x {\tilde v}\overline{{\tilde g}}_{3x}\, dx.
$$
Multiplying $(\ref{v4})_4$ by $ \overline{v}$ and integrating over $ (0,1)$ and using integration by parts we get
\begin{equation}\label{e36}
\Int_{0}^{1} a(x)|v_{x}|^{2} dx - \Int_{0}^{1} |{\tilde v}|^{2} dx-[a(x)v_x\overline{v}]_{0}^{1}
+\alpha\Int_{0}^{1}u {\overline v}\, dx= \Int_{0}^{1} {\tilde v}\overline{{\tilde g}_{3}}\, dx+\Int_{0}^{1}{\tilde g}_{4}\overline{v}dx.
\end{equation}
Multiplying (\ref{e36}) by $m_a/2$ and summing with (\ref{e35}) we get \\ \\
\begin{equation}
\matrix{\Int_{0}^{1} ((a(x)-xa'(x))+\frac{m_a}{2}a(x))|v_{x}|^{2}+(1-\frac{m_a}{2})|{\tilde v}(x)|^{2}]\, dx\hfill &\cr
+2\alpha\Re\Int_{0}^{1}xu {\overline v_x}\, dx+\alpha\Frac{m_a}{2}\Int_{0}^{1}u {\overline v}\, dx\hfill &\cr
=[x{\cal I}_{v}]_{0}^{1}+\frac{m_a}{2}[a(x)v_x\overline{v}]_{0}^{1}+R\hfill &\cr}
\label{ee22}
\end{equation}
with:
$$
R=R_1 + R_2
$$
and
$$
R_2= \frac{m_a}{2}\Int_{0}^{1} {\tilde v}\overline{{\tilde g}_{3}}\, dx+\frac{m_a}{2}\Int_{0}^{1}{\tilde g}_{4}\overline{v}dx.
$$
We have $[a(x)v_x\overline{v}]_{0}^{1}=0$ and $[x{\cal I}_{v}]_{0}^{1}=a(1)|v_x(1)|^2$.
Since $\|{\tilde v}\|_{L^2(0, 1)}, \|\sqrt{a(x)}v_x\|_{L^2(0, 1)}$ are bounded, we have from (\ref{ee22}):
\begin{equation}
a(1)|v_x(1)|^2\leq C.
\label{v6}
\end{equation}
By eliminating ${\tilde u}$ and ${\tilde v}$ from system (\ref{v4}) we obtain
\begin{equation}
\lambda^{2}u+(a(x)u_x)_x-\alpha v=f \hbox{ in } L^2(0, 1),
\label{v7}
\end{equation}
\begin{equation}
\lambda^{2} v+(a(x)v_x)_x-\alpha u=g \hbox{ in } L^2(0, 1),
\label{v8}
\end{equation}
where
\begin{equation}
\left\{\matrix{\|f\|_{L^2(0, 1)}=\left\|\Frac{g_2+i\lambda g_1}{\lambda^l}\right\|_{L^2(0, 1)}=\Frac{o(1)}{\lambda^{l-1}},\hfill & \cr
\|g\|_{L^2(0, 1)}=\left\|\Frac{g_4+i\lambda g_3}{\lambda^l}\right\|_{L^2(0, 1)}=\Frac{o(1)}{\lambda^{l-1}}.\hfill & \cr}\right.
\label{v9}
\end{equation}

Next we multiply (\ref{v7}) by $\overline{v}$ and (\ref{v8}) by $\overline{u}$, then add the resulting equations.
This yields
\begin{equation}
\matrix{\alpha\Int_{0}^{1}|v|^2\, dx=\alpha\Int_{0}^{1}|u|^2\, dx-\Re[av_x {\overline u}]_{0}^{1}
-\Re\Int_{0}^{1}i\lambda {\tilde u}{\overline v}\, dx+\Re\Int_{0}^{1}i\lambda {\tilde v}{\overline u}\, dx\hfill &\cr
\qquad\qquad+\Re\Int_{0}^{1}{\tilde g_2}{\overline v}\, dx-\Re\Int_{0}^{1}{\tilde g_4}{\overline u}\, dx.\hfill &\cr}
\label{v10}
\end{equation}
Then
\begin{equation}
\matrix{\alpha\Int_{0}^{1}|v|^2\, dx=\alpha\Int_{0}^{1}|u|^2\, dx-\Re[av_x {\overline u}]_{0}^{1}
+\Re\Int_{0}^{1}(i\lambda {\tilde g_1}+{\tilde g_2}){\overline v}\, dx-\Re\Int_{0}^{1}(i\lambda{\tilde g_3}+{\tilde g_4} ){\overline u}\, dx.\hfill &\cr}
\label{v80}
\end{equation}
Thus, applying Cauchy-Schwarz's and Young's inequalities, using (\ref{v5}) and (\ref{v6}) we obtain
\begin{equation}
\beta^2\Int_{0}^{1}|v|^2\, dx=\beta^2\Int_{0}^{1}|u|^2\, dx+\Frac{o(1)}{\lambda^{\frac{l-(1-\tau)}{2}-1}}.
\label{v11}
\end{equation}
\begin{lemma}
We have that
\begin{equation}
\matrix{\Int_{0}^{1} \left[\left((a(x)-xa'(x))+\frac{m_a}{2}a(x)\right)|u_{x}|^{2}+\left(1-\frac{m_a}{2}\right)|{\tilde u}(x)|^{2}\right]\, dx\hfill &\cr
+2\alpha\Re\Int_{0}^{1}xv {\overline u_x}\, dx+\alpha\Frac{m_a}{2}\Int_{0}^{1}v {\overline u}\, dx\hfill &\cr
=[x{\cal I}_{u}]_{0}^{1}+\frac{m_a}{2}[a(x)u_x\overline{u}]_{0}^{1}+R,\hfill &\cr}
\label{ee22wwn}
\end{equation}
where
$$
{\cal I}_{u}(x)= |\sqrt{a(x)}u_x(x)|^{2} +|{\tilde u}(x)|^{2}
$$
and
$$
R=2 \Re\Int_{0}^{1} x {\tilde g}_2  \overline{u}_x\, dx+2 \Re\Int_{0}^{1}x {\tilde u}\overline{{\tilde g}}_{1x}\, dx
+\frac{m_a}{2}\Int_{0}^{1} {\tilde u}\overline{{\tilde g}_{1}}\, dx+\frac{m_a}{2}\Int_{0}^{1}{\tilde g}_{2}\overline{u}dx.
$$
\label{ll15}
\end{lemma}
{\bf Proof}\\
To get (\ref{ee22wwn}), let us multiply the equation $(\ref{v4})_2$ by $x \overline{u}_x$ Integrating on $(0, 1)$ we obtain
$$
i\lambda\Int_{0}^{1} {\tilde u} x \overline{u}_x\, dx-\Int_{0}^{1} (a(x)u_x)_x x \overline{u}_x\, dx
+\alpha\Int_{0}^{1} v x \overline{u}_x\, dx
=\Int_{0}^{L}{\tilde g}_2 x \overline{u}_x\, dx
$$
or
$$
-\Int_{0}^{1} {\tilde u} x (\overline{i\lambda u_x})\, dx-\Int_{0}^{1} x (a(x)u_x)_x  \overline{u}_x\, dx
+\alpha\Int_{0}^{1} v x \overline{u}_x\, dx
= \Int_{0}^{1}{\tilde g}_2 x \overline{u}_x\, dx.
$$
Since $i\lambda u_x= {\tilde u}_x+f_{1x}$ taking the real part in the above equality results in
$$
\matrix{-\Frac{1}{2}\Int_{0}^{1} x \Frac{d}{dx}| {\tilde u}|^{2}\, dx+\Frac{1}{2}\Int_{0}^{1}x a(x)\Frac{d}{dx}|u_{x}|^{2}\, dx
-[x a(x) |u_{x}|^{2} ]_{0}^{1}+\Int_{0}^{1}  a(x)|u_{x}|^{2}\, dx
\hfill & \cr
+\alpha\Re\Int_{0}^{1} v x \overline{u}_x\, dx=   \Re\Int_{0}^{1 }{\tilde u} x\overline{{\tilde g}}_{1x}\,dx +  \Re\Int_{0}^{1}{\tilde g}_2 x \overline{u}_x\, dx .\hfill & \cr}
$$
Performing an integration by parts we get
\begin{equation}\label{e35z}
    \Int_{0}^{1}  [|\sqrt{a(x)} u_{x}|^{2}+|{\tilde u}(x)|^{2}]\, dx
-\Int_{0}^{1}  x a'(x)|u_{x}(x)|^{2}\, dx+2\alpha\Re\Int_{0}^{1} v x \overline{u}_x\, dx
=[x(|\sqrt{a(x)} u_{x}|^{2}+|{\tilde u}(x)|^{2})]_{0}^{1}+R_{1},
\end{equation}
where
$$
R_{1}=2 \Re\Int_{0}^{1} x {\tilde g}_2  \overline{u}_x\, dx+2 \Re\Int_{0}^{1}x {\tilde u}\overline{{\tilde g}}_{1x}\, dx.
$$
Multiplying $(\ref{v4})_2$ by $ \overline{u}$ and integrating over $ (0,1)$ and using integration by parts we get
\begin{equation}\label{e36z}
\Int_{0}^{1} a(x)|u_{x}|^{2} dx - \Int_{0}^{1} |{\tilde u}|^{2} dx-[a(x)u_x\overline{u}]_{0}^{1}
+\alpha\Int_{0}^{1} v \overline{u}\, dx= \Int_{0}^{1} {\tilde u}\overline{{\tilde g}_{1}}\, dx+\Int_{0}^{1}{\tilde g}_{2}\overline{u}dx.
\end{equation}
Multiplying (\ref{e36z}) by $m_a/2$ and summing with (\ref{e35z}) we get \\ \\
\begin{equation}
\matrix{\Int_{0}^{1} ((a(x)-xa'(x))+\frac{m_a}{2}a(x))|u_{x}|^{2}+(1-\frac{m_a}{2})|{\tilde u}(x)|^{2}]\, dx\hfill &\cr
+2\alpha\Re\Int_{0}^{1} v x \overline{u}_x\, dx+\alpha\Frac{m_a}{2}\Int_{0}^{1} v \overline{u}\, dx\hfill &\cr
=[x{\cal I}_{u}]_{0}^{1}+\frac{m_a}{2}[a(x)u_x\overline{u}]_{0}^{1}+R\hfill &\cr}
\label{ee22v}
\end{equation}
with:
$$
R=R_1 + R_2
$$
and
$$
R_2= \frac{m_a}{2}\Int_{0}^{1} {\tilde u}\overline{{\tilde g}_{1}}\, dx+\frac{m_a}{2}\Int_{0}^{1}{\tilde g}_{2}\overline{u}dx.
$$
\noindent
$\bullet${\bf Step 3}
We have $[x{\cal I}_{u}]_{0}^{1}=a(1)|u_x(1)|^2+|{\tilde u}(1)|^2$ and
$[a(x)u_x\overline{u}]_{0}^{1}=a(1)u_x(1){\overline u}(1)$.
By definition of $m_a$, we have
$$
(2-m_a)a\leq 2(a-xa')+m_a a.
$$
This, together with (\ref{ee22v}), gives
\begin{equation}
\Int_{0}^{1}(a(x)|u_x|^2+|{\tilde u}|^2)\, dx=\Frac{o(1)}{\lambda^{l-(1-\tau)}}.
\label{ee23zr}
\end{equation}
It follows from (\ref{v11}) and (\ref{ee23zr}) that
\begin{equation}
\Int_{0}^{1}|{\tilde v}|^2\, dx\rightarrow 0.
\label{ee23zz}
\end{equation}
Finally, from (\ref{e36}) and (\ref{ee23zz}), we obtain that
\begin{equation}
\Int_{0}^{1}a(x)|v_x|^2\, dx\rightarrow 0.
\label{ee23zzz}
\end{equation}
Since $\omega> 0$, we have
\begin{equation}
\|\varphi\|_{L^{2}(-\infty, \infty)}^2\leq \Frac{1}{\omega}\Int_{-\infty}^{+\infty}(\varsigma^{2}+\omega)|\varphi(\varsigma)|^{2}\, d\varsigma\rightarrow 0.
\label{ee23zzzz}
\end{equation}
Combining (\ref{ee23zr}), (\ref{ee23zz}),(\ref{ee23zzz}) and (\ref{ee23zzzz}), we obtain that
\begin{equation}
\|U\|_{{\cal H}}\rightarrow 0.
\label{ee23zzzzz}
\end{equation}
This is a contradiction with the assymption that $\|U\|_{{\cal H}}=1$.

Moreover the decay rate is optimal. In fact for the case $a(x)=x^{\gamma}, \gamma\in [0, 2[$,
the decay rate is consistent with
the asymptotic expansion of eigenvalues which shows a behavior of the real part like $k^{-(3-\tau)}$.

\hfill$\Box$\\

\end{sloppypar}

\begin{thebibliography}{99}


\bibitem{5} Z. Achouri, N. Amroun, A. Benaissa, {\em The Euler-Bernoulli beam
equation with boundary dissipation of fractional derivative type. Math\/,}
Method. Appl. Sci. {\bf 40} (2017)-11, 3837-3854.

\bibitem{akil} M. Akil, M. Ghader and A. Wehbe, {\em The influence of the coefficients of a system of wave equations
coupled by velocities on its stabilization\/,} SeMA J. {\bf 78} (2021)-3, 287-333.

\bibitem{alabau1} F. Alabau-Boussouira, P.Cannarsa \& G. Leugering, {\em Control and stabilization of degenerate
wave equations\/,} {SIAM J.Controle Optim}, {\bf 555}, (2017)-3, 1-36.


\bibitem{arba} W. Arendt and C. J. K. Batty, {\em Tauberian theorems and stability of one-parameter semigroups\/,}
Trans. Am. Math. Soc., {\bf 306} (1988), pp. 837-852.
\bibitem{beai} A. Benaissa, C. Aichi, {\em Energy decay for a degenerate wave equation under
fractional derivative controls\/,} Filomat {\bf 32} (2018)-17, 6045-6072.

\bibitem{boto} A. Borichev and Y. Tomilov, {\em Optimal polynomial decay of functions and operator semigroups\/,}
{Math. Ann.} {\bf 347} (2010)-2, 455-478.

\bibitem{canna.3} P. Cannarsa, P. Martinez and J. Vancostenoble, {\em Carleman estimates for a class
of degenerate parabolic operators\/,} {SIAM J. Control Optim.,} {\bf 47}, (2008)-1, 1-19.
(electronic), 2006.

\bibitem{choi} J. U. Choi and R. C. Maccamy, {\em Fractional order Volterra equations with applications to elasticity\/,}
J. Math. Anal. Appl., {\bf 139} (1989), 448-464.

\bibitem{foto} M. Fotouhi, L. Salimi, {\em Null controllability of degenerate/singular parabolic equations\/,}
J. Dyn. Control Syst.  {\bf 18}  (2012)-4, 573-602.

\bibitem{13} F. Huang, {\em Characteristic conditions for exponential stability
of linear dynamical systems in Hilbert spaces\/,} Ann. Differ. Equ., {\bf 1} (1985),
43-55.
\bibitem{1} M. Kerdache,M. Kesri, A. Benaissa, {\em Fractional boundary
stabilization for a coupled system of wave equations\/,} Ann. Univ. Ferrara
Sez. VII Sci. Mat. {\bf 67} (2021)-1, 121-148.
\bibitem{weko} M. Koumaiha, {Analyse num\'erique pour les \'equations de Hamilton-Jacobi sur r\'eseaux
et contr\^{o}labilit\'e / stabilit\'e indirecte d'un syst\`eme d'\'equations des ondes 1D\/,}
PhD thesis, Université Paris est, 2017.
\bibitem{16} V. Komornik. Exact Controllability and Stabilization : The
Multiplier Method. Wiley-Masson Series Research in Applied Mathematics.
Wiley, 1995.
\bibitem{lebed} N. N. Lebedev, Special Functions and their Applications, Dover
Publications, New York, (1972).
\bibitem{lira} Z. Liu, B. Rao, {\em Frequency domain approach for the polynomial stability of a system
of partially damped wave equations\/,} J. Math. Anal. Appl. {\bf 335} (2007)-2, 860-881.
\bibitem{4} B. Mbodje, {\em Wave energy decay under fractional derivative
controls\/,} IMA J. Math. Contr. Inf. {\bf 23} (2006), 237-257.
\bibitem{9} A. Pazy, Semigroups of Linear Operators and Applications to
Partial Differential Equations. Springer Verlag, New York (1983).
\bibitem{12} J. Pruss, {\em On the spectrum of C0-semigroups\/,} Transactions of the
American Mathematical Society, {\bf 284} (1984)-2, 847-857.
\bibitem{zerk} H. Zerkouk, C. Aichi, A. Benaissa, {\em  On the stability of a degenerate wave equation
under fractional feedbacks acting on the degenerate boundary\/,} J. Dyn. Control Syst. {\bf 28} (2022)-3, 601-633.




















\end{thebibliography}
\end{document}